\documentclass[11pt]{amsart}

\usepackage[usenames]{color}
\usepackage{amsmath,amsfonts,amssymb,epsfig}
\usepackage{latexsym}
\usepackage{amsthm}
\usepackage{hyperref}
\usepackage{cancel}
\usepackage{mathrsfs}
\usepackage{palatino}

\newtheorem{theorem}{Theorem}
\newtheorem{definition}[theorem]{Definition}
\newtheorem{remark}[theorem]{Remark}
\newtheorem{corollary}[theorem]{Corollary}
\newtheorem{proposition}[theorem]{Proposition}
\newtheorem{lemma}[theorem]{Lemma}

\newtheorem{question}[theorem]{Question}
\newtheorem{example}{Example}

\newcommand{\Mg}{D_g}
\newcommand{\mg}{d_g}

\newcommand{\no}{\noindent}

\newcommand{\conv}{\text{\rm conv}}

\newcommand{\bex}{\begin{example}\em}
\newcommand{\eex}{\end{example}}

\def\R{\mathbb{R}}

\def\dfn#1{{\em #1}}

\def\SU2{\mathrm{SU}(2)}
\def\half{\frac{1}{2}}

\def\inj{\mathrm{inj}}
\def\sec{\mathrm{sec}}
\def\tr{\mathop{\mathrm{tr}}}
\def\II{I\!I}

\newcommand\g[1]{\left\langle #1 \right\rangle}
\def\F{\mathscr{F}}  
\def\xistd{\xi_\text{std}} 

\renewcommand\div{\mathrm{div}}
\def\vol{\mathrm{vol}}

\def\S{\mathrm{S}}
\def\diam{\mathrm{diam}}

\def\cprime{$'$} 

\graphicspath{{./pict/}}

\numberwithin{equation}{section} 
\numberwithin{theorem}{section}

\title[Tightness in contact metric 3--manifolds]{Tightness in contact metric
3--manifolds} 

\author{John B. Etnyre} 

\address{School of Mathematics,
Georgia Institute of Technology,
Atlanta, GA 30332} \email{etnyre@math.gatech.edu}
\urladdr{http://www.math.gatech.edu/\~{}etnyre}

\author{Rafal Komendarczyk}

\address{Department of Mathematics,
Tulane University,
New Orleans, LA 70118 } \email{rako@tulane.edu}
\urladdr{http://www.math.tulane.edu/\~{}rako}

\author{Patrick Massot}
\address{Universit\'e Paris Sud 11, 91430 Orsay, France}
 \email{patrick.massot@math.u-psud.fr}
\urladdr{http://www.math.u-psud.fr/\~{}pmassot/}

\begin{document}
\begin{abstract}
This paper begins the study of relations between Riemannian geometry and
global properties of contact structures on 3--manifolds.
In particular we prove  an analog of the sphere theorem from Riemannian
geometry in the setting of  contact geometry. 
Specifically,  if a given three dimensional contact manifold $(M,\xi)$ admits a
complete compatible Riemannian metric of positive $4/9$-pinched curvature then
the underlying contact structure $\xi$ is tight; in particular, the contact
structure pulled back to the universal cover is the standard contact structure
on $S^3$.  We also describe geometric conditions in dimension three for $\xi$
to be universally tight in the nonpositive curvature setting.  
\end{abstract} 

\maketitle

\vspace{-1cm}

\section{Introduction}

A contact structure on 
a 3--manifold is either overtwisted or tight according as it contains, or does 
not contain, an embedded disk which is tangent to the contact planes 
along its boundary \cite{Eliashberg_vrille}. 
The study of the tight/overtwisted dichotomy and its implications for
low-dimensional topology is the primary
focus of contact topology, one of the newest branches of contact geometry
\cite{Bennequin83,Eliashberg92a}. A driving question in the area concerns the
existence of tight contact structures on a given manifold.  Though this has
been much studied it is still very difficult to distinguish tight from
overtwisted contact structures. In this paper we present several results
relating tightness of contact structures to curvature constraints on metrics
suitably adapted to the contact structure. In particular, we also introduce 
new ways for a Riemannian metric and a contact structure to be compatible
(though our notion of ``weakly compatible'' discussed below has had 
precursors in the literature \cite{EtnyreGhrist00, Rafal08}, it does not seem 
to have been formalized until
now).

A natural reference point here is the well known sphere theorem in Riemannian
geometry \cite{Rauch51, Berger60, Klingenberg61}, which is one of the
fundamental results showing how geometry can control the topology of the domain.
Recall, the sphere theorem states that every simply connected $n$-manifold
which admits $1/4$-pinched positive sectional curvature is homeomeorphic to the
$n$-sphere. In a similar vein we ask if a suitably chosen class of Riemannian
metrics can control the topology of an underlying contact structure. Slightly
extending a classical definition \cite{ChernHamilton85}, we say that a metric
is compatible with a contact structure if the corresponding unit contact form
$\alpha$ satisfies $*d\alpha = \theta' \alpha$ where $*$ is the Hodge star
operator induced by the metric and $\theta'$ is some constant.  We establish the
following theorem (in dimension three), which we refer to as the {\em contact
sphere theorem}.
%
\begin{theorem}[Contact sphere theorem]
\label{CST-thm}
Let $(M,\xi)$ be a closed contact 3--manifold and $g$ a complete Riemannian
metric compatible with $\xi.$ If there is a constant $K_{max}>0$ such that the
sectional curvatures of $g$ satisfy 
\[
0<\frac 49 K_{max}<\sec(g)\leq K_{max},
\]
then the universal cover of $M$ is diffeomorphic to the 3--sphere by a
diffeomorphism taking the lift of $\xi$ to the standard contact structure on
the 3--sphere. 
\end{theorem}
Recall the standard contact structure $\xistd$ on the  3--sphere $S^3$ is the
one induced as the complex tangencies to the unit sphere in $\mathbb{C}^2.$
This contact structure can alternately be described as the orthogonal
complement, in the round metric on $S^3,$ to the Hopf fibration. Eliashberg
proved in \cite{Eliashberg92a} that $\xistd$ is the unique (up to
contactomorphism)  tight contact structure on $S^3$, and what we actually prove
in the above theorem is tightness of the contact structure. So this theorem
really gives Riemannian geometric conditions that imply tightness.

Turning to general 3--manifolds, we recall that there are
manifolds that do not admit tight contact structures \cite{EtnyreHonda01a} and
we know precisely which Seifert manifolds admit them
\cite{LiscaStipsicz09}. Since we also know manifolds containing homologically
essential surfaces admit tight contact
structures \cite[Corollary~3.2.11]{EliashbergThurston98} we are left with the
central question: Which hyperbolic manifold (specifically, homology spheres)
can support a tight contact structure?

Many hyperbolic manifolds are known to have tight contact structures but their
construction and the proof of their tightness has nothing to do with the
hyperbolic metric. While we do not directly address this question here, it does
provide a strong motivation for better understanding the relation between
contact topology and Riemannian geometry. We also state some results that
provide a hint as to how to approach this and similar questions.  For example,
as a corollary of Theorem \ref{thm:weak-compatible} below, we obtain a
geometric criterion for a contact structure on a negatively curved manifold to
be universally tight, see Theorem~\ref{thm:hyp-criterion}. (A contact structure
is universally tight if its pull back to the universal cover of the manifold is
tight, this is stronger than tightness.)

From a slightly different perspective, one might also want to go beyond the tight
vs overtwisted dichotomy and sort the class of tight contact structures by
finding privileged subclasses. In addition to the notion of universal
tightness explained above, a classical class is that of symplectically fillable
contact structures. In
this paper, we seek classes interacting nicely with curvature in Riemannian
geometry. Note that such interactions will automatically be inherited by
covering spaces, contrasting with fillability properties.

Our approach to the contact sphere theorem and other global questions of
tightness is through a quantitative versions of Darboux theorem. Recall that
Darboux theorem in contact geometry guaranties that each point in a contact
manifold has a neighborhood which is standard, i.e. embeds inside the standard
contact structure on $\R^3$.  According to Bennequin's theorem this
neighborhood is then tight. In a Riemannian setting we can ask for a
quantitative version guarantying that balls up to a certain radius are standard
and, in particular, tight. 
Let $g$ be any Riemannian metric and $\xi$ any contact
structure on $M.$ We define the \dfn{tightness radius of $\xi$ and $g$ at a
point $p$} to be
\begin{align*}
 \tau_p(M,\xi)  =\sup\{r\,|\,& \text{the contact structure on the geodesic } \\
 			&\text {ball $B_p(r)$ at $p$ of radius $r$ is tight}\},
\end{align*}
and the global \dfn{tightness radius}
\[
 \tau(M,\xi)  =\inf_{p\in M} \tau_p(M,\xi).
\]
Of course, if we do not assume any compatibility condition between the metric
and the contact structure then we cannot estimate the tightness radius. 
We first concentrate on what happens with the compatibility definition recalled
above.
The tightness radius is, by definition, always less than the injectivity
radius but one could ask if, for compatible metrics, they always coincide. 
This would explain the following surprising result (which is an important
ingredient of the proof of Theorem \ref{CST-thm}). 
\begin{theorem}\label{thm:onspheres}
Let $(M, \xi)$ be a contact 3--manifold and $g$ a complete Riemannian metric
that is compatible with $\xi.$
For a fixed point $p\in M$ let $\tau_p= \tau_p(M,\xi)$ and suppose that
$\tau_p < \inj_p(g).$ Then for all radii $r$ with $\tau_p \leq r < \inj_p(g)$,
the geodesic sphere $S_p(r)$ contains an overtwisted disk.
\end{theorem}

\no
Recall that, a priori, overtwisted disks can have a very complicated geometry
and this is what makes it hard to prove tightness of contact structures. We
find this theorem somewhat surprising as it says that when a metric is
compatible with a contact structure then as soon as a geodesic ball is large
enough to be overtwisted one sees the overtwisted disk in a specific place,
namely the boundary of the ball. Thus making it easy to determine when such a
ball is tight (using Bennequin's theorem).

Despite this surprising result, we have numerical experiments, see Section \ref{S:overt-num}, 
which strongly
suggest that the tightness radius can indeed be less than the injectivity
radius for compatible metrics, so we search for geometrical quantities
controlling the tightness radius. To this end we recall that given a Riemannian
metric $g$ on $M$ the \dfn{convexity radius} of $g$ is defined to be 
\[
\begin{split}
\conv(g)=\sup \{ r \, |\,& r< \inj(g) \text{ and the geodesic balls of radius $r$} \\ 
&\text{are weakly geodesically convex} \},
\end{split}
\]
where $\inj(g)$ is the injectivity radius of $(M,g).$ For a more complete
discussion see Subsection~\ref{ss:riemannianconvexity}.

\begin{theorem}\label{thm:compatible}
Let $(M, \xi)$ be a contact 3--manifold and $g$ a complete Riemannian metric 
that is compatible with $\xi.$ Then,
\begin{equation}\label{tightest}
 \tau(M,\xi)\geq \conv(g).
\end{equation}
In particular, if $\sec(g)\leq K$, for $K > 0$,
then
\[
 \tau(M,\xi)\geq \min\{\inj(g),\frac{\pi}{2\sqrt{K}}\}
\]
and $\tau(M,\xi) = \inj(g),$ if $g$ has non-positive curvature.
\end{theorem}

We note that if $M$ is a compact manifold then one may easily show that a lower
bound for $\tau(M,\xi)$ exists. To see this note that  $M$ may be covered with
Darboux balls (which are tight). Then the Lebesgue number for this open cover
provides the desired lower bound. Of course this bound exists for any metric and
one has virtually no control over it. Theorem~\ref{thm:compatible} shows that if the
metric is compatible with the contact structure then one does not need compactness
and one can estimate $\tau(M,\xi)$ below in terms of curvature and injectivity
radius.  In particular, our theorem shows when $M$ is noncompact with bounded
curvature and injectivity radius, the tightness radius cannot shrink to zero at
infinity.

The above theorem is based on the comparison of Riemannian convexity and
almost-complex convexity in symplectizations of contact manifolds. Its proof
uses holomorphic curves techniques due to Gromov and Hofer. In the seminal paper
\cite{Hofer93}, Hofer proved that overtwisted disks guaranty the existence of
closed Reeb orbits. Here  almost-complex convexity controlled by
Riemannian convexity allows us to use this argument in balls where we know there
are no closed Reeb orbits. (In
Section~\ref{S:xiS}, we will see a surprising link with a third type of
convexity --- Giroux's convexity of surfaces in contact manifolds.)

We notice that our bounds on the tightness radius are especially effective in
the case of non-positive curvature. 
\begin{corollary}\label{cor:ut}
Let $(M,\xi)$ be a contact manifold and $g$ a complete Riemannian metric
compatible with $\xi$ having non-positive sectional curvature. Then $\xi$ is
universally tight and hence tight.
\end{corollary}

One should point out that the study of compatible metrics is useful in fluid
mechanics, plasma physics and other subjects, see for example
\cite{EtnyreGhrist00}.  In addition, it has produced a great many questions from
the Riemannian geometry perspective, see for example \cite{Blair02}. 
It appears however, that so far there has been little work connecting properties
of compatible metrics with much studied global properties of the contact
structure in dimension three, such as tightness (though a exception to this is
\cite{Rafal08} that provides a lower bound for the volume of overtwisted Seifert
fibered manifolds under $S^1$-symmetry conditions).

The class of compatible metrics is very natural but it is fairly restrictive.
In relation to the question of existence of tight contact structure on
hyperbolic manifolds, we note right away that a hyperbolic metric cannot be
compatible with a contact structure on a closed manifold \cite[p.\ 99]{Blair02}.
Furthermore, Blair conjectures that if a metric is compatible with a contact
structure on a closed manifold and has non-positive curvature then it is flat.
So the above corollary may be of very limited impact\footnote{However it does
not seem obvious that any contact structure which is compatible with a
flat metric is tight, especially since there is no classification of foliations
of $\R^3$ by lines, contrasting with the situation in $\S^3$.}.  This justifies
the introduction of a more general class of metrics that can include the
hyperbolic ones.  Another motivation for extending the notion of compatibility
comes from the theory of curl eigenfields. A curl eigenfield on a Riemannian
3--manifold is a 1--form $\alpha$ satisfying $*d\alpha = \theta' \alpha$
with $\theta'$ constant but $\|\alpha\|$ can vary.  (We note that this equation
is dual to the normal curl eigenfield equation for vector fields.) See
\cite{EtnyreGhrist00} for some applications of this concept.

We say that a Riemannian metric and a contact structure $\xi$ are
\dfn{weakly compatible} if there exist a Reeb vector field for $\xi$ which is
perpendicular to $\xi$. 
We will show in Section~\ref{S:compatible_metrics} that this condition can be
equivalently stated as there exists a contact form $\alpha$ such that
\begin{equation}
\label{eqn:compatible}
*d\alpha = \theta'\alpha,
\end{equation}
where $\theta'$ is a positive function (which we will see measures the rotation
speed of the contact planes).
This equivalent definition enables us to see that the class of weakly compatible
metrics is an extension of compatible metrics.

This class of metrics 
includes all the non-singular curl eigenfields and Beltrami fields. In addition, 
as shown in Section~\ref{S:examples} 
it allows for hyperbolic metrics.  It is also
stable both under conformal changes and under the modifications used by Krouglov
in \cite{Krouglov08}, see Remark~\ref{rem:Krouglov} below.

We will use several measures of how far a weakly compatible metric is from being
compatible. First the rotation speed $\theta'$ and the norm $\rho$
of the special Reeb vector field $R$ entering in the definition are both
constant in the compatible case so their gradient are such measures. We will
also use the mean curvature $H$ of the contact plane field. It vanishes for
compatible metrics and its definition is recalled in
Section~\ref{S:compatible_metrics}. Finally we shall also use the normalized
Reeb vector field $n = R/\|R\|$, which also happens to be a unit normal vector 
field to the contact planes. In the compatible case it is a geodesic vector
field so we will consider $\nabla_n n$.

We will prove in Section~\ref{S:compatible_metrics} that the following two
combinations of these measures give the same vector field:
\begin{equation}\label{eq:M_g}
\Mg := \nabla_n n + 2H n = (\nabla \ln\theta')^\perp - \nabla \ln\rho
\end{equation}
where $v^\perp$ is the component of $v$ perpendicular to $\xi.$
We introduce:
\begin{equation*}\label{eq:m_g}
 \mg = \max_M \|\Mg\| 
\end{equation*}

Note that $\mg$ is finite whenever $M$ is compact and vanishes for compatible
metrics (all terms in $\Mg$ vanish in this case).

To extend our main theorem to weakly compatible metrics we also introduce the
following notation: let  $K\geq 0$ and 
$\sec(g)\leq \pm K$, define 
 \begin{equation}\label{cotsdef}
 \begin{split}
  \text{\rm ct}_K(r) & =\begin{cases}
  \sqrt{K}\cot(\sqrt{K} r)\, , & \qquad \text{for $\sec(g)\leq K$,\ \ $r\leq \min\{\inj(g),\frac{\pi}{2\sqrt{K}}\}$},\\
  \frac{1}{r}, & \qquad \text{for $\sec(g)\leq 0$},\\
  \sqrt{K}\coth(\sqrt{K} r), &\qquad \text{for $\sec(g)\leq -K$.}
  \end{cases}\\
 \end{split}
\end{equation}
Here, of course, in the first case we assume $\sec(g)$ is positive
somewhere and in the second case that it is 0 somewhere. Also to simplify our
notations we will often write $\text{\rm ct}_{K}$ instead of $\text{\rm
ct}_{-K}$ understanding that we mean the latter in the negative curvature
setting.  We may now state our result for weakly compatible metrics as follows. 

\begin{theorem}
\label{thm:weak-compatible}
 Let $(M,\xi)$ be a contact 3-manifold (not necessarily closed) that is weakly
 compatible with a Riemannian metric $g$. Whenever $\mg < \infty$ the tightness radius
 admits the following lower bound
\[
 \tau(M,\xi)\geq \min\left\{\text{\rm ct}^{-1}_{K}\left(\mg\right), \inj(g)\right\}.
\]
\end{theorem}

\begin{remark}
{\rm
While the above theorem provides a bound  on the tightness radius in a weakly 
compatible metric, we note that it is not sufficient to prove a version of the the 
Contact Sphere Theorem~\ref{CST-thm} in this setting; however, it is reasonable 
to hope that such a theorem holds for weakly compatible metrics, but possibly 
with weaker bounds. Evidence for this, as well as a discussion of other possible 
strengthenings of Theorem~\ref{CST-thm}, is discussed in 
Section~\ref{sec:extensions}.
}
\end{remark}

\begin{remark}
{\rm
\label{rem:Krouglov}
The above theorems can be applied only when we have control over the sectional
curvature of all plane fields and not only the sectional curvature of contact
planes $\xi$.  This is natural in view of the following very slight sharpening of
a result of Krouglov saying that the latter curvature is very 
flexible. To get this version, start with a compatible metric and observe that
Krouglov's modifications do not destroy weak compatibility, although they destroy
compatibility (this is another reason to use weakly compatible metrics).

\begin{theorem}[Krouglov 2008, \cite{Krouglov08}]
Given a cooriented contact structure $\xi$ on a closed 3--manifold $M$ and any
strictly negative function $f,$ there is a weakly compatible metric on $M$ such
that the sectional curvatures of $\xi$ are given by $f$. Moreover, if the Euler
class of $\xi$ is zero then any function $f$ may be realized. \qed
\end{theorem}
}
\end{remark}

\no Observe that since $\text{ct}_K(r)\to\infty$ as $r\to 0$ the bound in
Theorem \ref{thm:weak-compatible} is always
nonzero. We also notice that if $\alpha$ is actually compatible with $g$ then
$\mg = 0$ and thus $\text{\rm ct}^{-1}_{K}(\mg)$ can be taken to be $+\infty$. A
similar situation occurs when $\sec(g)\leq -K$, then $\text{\rm
ct}^{-1}_{K}(r)$ is ill defined for $r\in [-\sqrt{K},\sqrt{K}]$ and we may
assume $\text{\rm ct}^{-1}_{K}(\mg)$ to be $+\infty$ as well. Recall that for
such manifolds, the universal cover is exhausted by geodesic balls. Since an
overtwisted disk has to be contained in a compact part of the universal cover,
we get the following corollary that we state as a theorem due its potential
relevance to the problem of the existence of tight contact structures on 
hyperbolic manifolds.
\begin{theorem}\label{thm:hyp-criterion}
 Let $(M,\xi)$ be a contact 3-manifold (not necessarily closed) that is weakly
 compatible with a complete Riemannian metric $g$ of nonpositive sectional
 curvature.  If
 \begin{equation}\label{eq:m_g-K}
   \sec(g) \leq -\mg^2
 \end{equation}
at all points then the contact structure $\xi$ is universally tight. 
\end{theorem}

One remarkable property of compatible metrics is that Reeb orbits are geodesics
(see Corollary~\ref{Rderivative})
and we use this in our study of compatible metrics. However, this is precisely
what rules out closed hyperbolic manifolds in dimension 3: these manifolds cannot have any
geodesic vector field \cite{Zeghib}. But many hyperbolic manifolds have
quasi-geodesic vector fields, see \cite{Calegari} for a recent account. These
vector fields also cannot have any contractible Reeb orbits. So if a closed
hyperbolic manifold has a quasi-geodesic Reeb field then the corresponding
contact structure is universally tight. This observation does not explicitly use
any easily defined compatibility between a metric and contact structure.
However we can use an easy differential geometric criterion for
quasi-geodesicity to get the following theorem which can then be compared to
Theorem~\ref{thm:hyp-criterion}.
\begin{theorem}\label{thm:vwc}
Let $(M,\xi)$ be a closed contact manifold. Suppose $M$ admits a metric $g$ such
that the sectional curvature of $g$ is bounded above by $-K$ for some constant
$K \geq 0$ and there is a Reeb vector field $R$ for $\xi$ such that the normalized
Reeb field  $N=R/\|R\|$ satisfies 
\[
\|\nabla_N N\| \leq \sqrt K.
\]
Then the universal cover of $M$ is tight.

\end{theorem}
We note that one can think of the condition $\|\nabla_N N\|\leq \sqrt K$ as some
type of compatibility between $g$ and $\xi$. We also note that, while this
theorem is stronger than Theorem~\ref{thm:hyp-criterion} when they both apply,
it does require that we are working with a closed manifold. We would lastly like
to point out that earlier we used $n$ for the normalized Reeb vector field while here
we used $N$ for that purpose. We will always use $n$ to denote a unit normal 
vector field to the contact planes (which, in a weakly compatible metric, the 
normalized Reeb vector field always is) and use $N$ if the normalized field does 
not have to be normal.

\medskip

\no {\em Outline:} The rest of the paper is structured as follows.
Section~\ref{S:compatible_metrics} defines the various notions of compatibility
between metrics and contact structures which will be used in this paper and
proves some formulas useful for our convexity comparison results.
Section~\ref{S:convexity_comparison} 
contains our results comparing Riemannian convexity and almost complex
convexity (Propositions~\ref{prop:levi-form_weak} and \ref{prop:levi-form}) and
the proofs of Theorems \ref{thm:compatible} and \ref{thm:weak-compatible} 
as well as
their corollaries. Section~\ref{S:quasi-geod} proves Theorem~\ref{thm:vwc}.
Section~\ref{S:xiS} centers around characteristic foliations of surfaces in
contact 3--manifolds which are used in the proof of
Theorem~\ref{thm:onspheres}.
Section~\ref{S:CST} contains the proof of the contact sphere theorem and
further discussion around it. 
Section~\ref{S:examples} describe examples where we can apply our estimates of
the tightness radius.

We end the introduction by noting that there are generalizations of some of the 
theorems presented above to the higher dimensional setting. These appear in 
\cite{EKM-higher} along with other methods of proof for some weaker forms of the 
Darboux theorem with estimates considered here. 

\medskip
\no {\emph{Acknowledgements}:}  The authors wish to thank Vladimir Krouglov for
useful correspondence and Fr\'ed\'eric Le Roux for his help with Lemma
\ref{technicallemma}. We also thank the referee for useful comments that improved 
the exposition of the paper. The second author is grateful to Chris Croke for
enjoyable lunch meetings, and many helpful suggestions to the project. The first
author was partially supported by the NSF Grant DMS-0804820. The second author
was partially supported by DARPA grant \#FA9550-08-1-0386.

\section{Compatible metrics}
\label{S:compatible_metrics}

In this section we discuss various notions of compatibility between a metric and
a contact structure. Our notions of compatibility are more general than the usual
notions of compatibility, but seem to be more natural from a contact geometry
perspective. 
Several of the results in this section come from, or were inspired by
(previously unpublished) parts of \cite{Patrick_thesis}.

Throughout this section and the rest of the paper we will interchangeably use
the notation $g(u,v)$ and $\g{u,v}$ to denote the Riemannian metric evaluated on
the vectors $u$ and $v.$

\subsection{Instantaneous rotation}

We begin by observing how to detect if a plane field is contact using any
Riemannian metric on a 3--manifold. To this end let $\xi$ be an arbitrary
oriented plane field on an oriented 3--manifold $M$. Recall that the Frobenius
integrability criterion says that $\xi$ is integrable if and only if the flow of
a (local) non-zero vector field tangent to $\xi$ preserves $\xi.$ Given a
Riemannian metric $g$ on $M$ we can (locally) choose an oriented orthonormal
moving frame $u, v, n$ where $u,v$ is an oriented basis for $\xi$ and $n$ is a
unit normal vector to $\xi.$ Denote by $\phi_t$ the flow of $u$ and define
$\theta(t)$ to be the angle between $(\phi_{-t})_*v$ and $n.$ That is
\[
\cos \theta(t)= \frac{g((\phi_{-t})_*(v), n)}{\| (\phi_{-t})_*(v) \|}.
\]
One may compute that $\theta'(0)=-g([u,v], n)$. Setting $\alpha(\cdot) = g(n, \cdot)$,
one can then characterize $\theta'(0)$ by 
$\alpha \wedge d\alpha = \theta'(0) \vol_g$. In particular, $\theta'(0)$
depends on $g$ and $\xi$ but not on $u$,~$v$. 

We denote the function $\theta'(0)$ by $\theta'$ and call it the \dfn{rotation}
or \dfn{instantaneous rotation} of $\xi$ with respect to $g.$ Notice that the
Frobenious condition implies that $\xi$ is a (positive) contact structure if and
only if $\theta'>0.$

\subsection{Second fundamental form}
In analogy with foliations we now introduce the second fundamental form of a
general plane field. This notion goes back at least as far as \cite{Reinhart}.
The second fundamental form of $\xi$ is the quadratic form on $\xi$
defined as follows: for vectors $u$ and $v$ in $\xi_p=T_pM\cap \xi,$ 
\begin{equation}\label{definitionII}
\II(u, v) = \frac{1}{2}\g{\nabla_u v+\nabla_v u, n},
\end{equation}
where $n$ is the oriented unit normal to $\xi.$ (We note that $u$ and $v$ will
need to be extended to vector fields tangent to $\xi$ in a neighborhood of $p$
to compute $\II(u,v),$ but the value of $\II(u,v)$ is independent of this
extension {\em i.e.\ }$\II$ is tensorial.)

We note the following geometric interpretation of $\II$ inspired by
\cite{Giroux_osculating} and proved in \cite{Patrick_thesis}.
For any point $p\in M$  and a unit vector $v\in \xi$ let $P_v$ denote the plane
spanned by $v$ and the oriented unit normal to $\xi.$ There is a neighborhood
$N$ of the origin in $T_pM$ such that the exponential map pulls $\xi$ back to a
plane field $(\exp\! |_N)^*\xi$ that is transverse to $P_v\cap N.$ This plane
field induces a foliation on $P_v\cap N$ and $\II(v,v)$ is the curvature of the
leaf of this foliation through the origin (measured by the flat metric on $T_pM$
given by $g_p$). 

We have the following curvatures derived from $\II.$ The \dfn{extrinsic
curvature}, $K_e$, of $\xi$ is the determinant of the quadratic form $\II$ with
respect to $g.$ The \dfn{mean curvature}, $H$, of $\xi$ is half the trace of
$\II$ (that is, it is the mean of the eigenvalues of $\II$).  It is clear that
when $\xi$ is an integrable plane field, then $\II$ is the standard second
fundamental form of the leaves of the foliation associated to $\xi$ and the quantities $K_e$
and $H$ are their classical curvatures.

\subsection{Weak compatibility}
We now want to introduce a notion of weak compatibility between a contact
structure and a metric.
For the sake of comparisons with other kinds of compatibility and for later
computations, we describe several characterizations of weak compatibility.

\begin{proposition}
	\label{prop:charac_wc}
Let $\alpha$ be a positive contact form on a 3--manifold $M$ and $g$ be a
metric on $M$. We set $\rho = \|R_\alpha\|$.
The following properties are equivalent:
\begin{enumerate}
\item
The Reeb vector field $R_\alpha$ is orthogonal to $\xi$.

\item
There is some positive function $\theta'$ such that
\[
*d\alpha=\theta'\alpha,
\]
where $*$ is the Hodge star operator associated to $g$. 
(If this is true then $\theta'$ is the instantaneous rotation
of $\ker \alpha$ with respect to $g.$)

\item 
There is an almost complex structure $J$ on $\xi$
that is compatible with $d\alpha|_{\xi},$ (that is
$d\alpha(Ju$, $Jv)$ $=d\alpha(u,v)$ for all $u,v\in \xi$ and $d\alpha(v,Jv)>0$ for
all $v\not=0$ in $\xi$) such that, using the extension $\phi$ of $J$ to $TM$ sending
$R_\alpha$ to 0, we have:
\[
g(u,v)=\frac{\rho}{\theta'} d\alpha(u,\phi (v)) + \rho^2 \alpha(u)\alpha(v).
\]
where $\theta'$ is the instantaneous rotation of $\ker\alpha$.
\end{enumerate}
\end{proposition}

\begin{definition}
{\em
We say a contact form and a Riemannian metric $g$ on a 3--manifold $M$ are
\dfn{weakly compatible} if the conditions of the preceding proposition are
satisfied. We say that a contact structure $\xi$ and a metric $g$ are weakly
compatible if there is a contact form $\alpha$ for $\xi$ that is weakly
compatible with $g$.}
\end{definition}

The first condition of the proposition is the most concise so it can be taken
as \emph{the} definition of weak compatibility. The second condition clearly
relates this notion to \cite{ChernHamilton85} which asks in addition that
$\rho = 1$ and $\theta' = 2$ everywhere. An intermediate notion between this
and our weak compatibility is when $\theta'$ is any constant but
$\rho$ can vary. Then one says that $\alpha$ is a \dfn{non-singular curl
eigenfield}. One also says that $\alpha$ is a \dfn{non-singular Beltrami field}
if $\theta'$ is arbitrary and $\delta \alpha=0,$ where $\delta$ is the formal
adjoint of the de Rham differential.  (The condition $\delta \alpha=0$  is
equivalent to the metric dual of $\alpha$ being divergence free.) These two
notions of compatibility, which can be generalized to arbitrary plane fields,
have been made use of in the study of ideal fluid dynamics, see for 
example \cite{EtnyreGhrist00}. 

The third condition relates to the alternative presentation of the notion of
\cite{ChernHamilton85} which can be found in \cite{Blair02}.
(Also notice that in much of the literature about contact metric geometry, in
particular \cite{Blair02}, contact structures are taken to be negative contact
structures on 3--manifolds --- that is, they induce the opposite orientation to
the one given on the manifold.  This introduces various signs that differ from
ours due to the fact that we take positive contact structures.) 

The reason for weakening the definition of compatibility is not only to widen
the class of metrics for which meaningful theorems can be proved but also to
allows for homothety as well as other alterations of compatible metrics. More
specifically, the set of metrics that are weakly compatible with a given $\xi$
is closed under conformal transformations as well as changes of the metric
restricted to $\xi$ and $\xi^\perp.$ 

We will only explain the equivalence of (1) and (2) since
we will not crucially use (3).

\begin{proof}
We first prove that (1) implies (2).
Let $v,w$ be a (local) oriented orthonormal frame for $\xi,$ $\rho=\|R_\alpha\|$
and $n=\frac{1}{\rho}R_\alpha.$ So $v,w, n$ is an oriented orthonormal frame for
$TM.$ Let $v^\ast, w^\ast,n^\ast$ be the dual basis. Clearly 
$\rho\, \alpha= n^\ast$ since the two 1--forms agree on our chosen basis. Since
$n$ is parallel to $R_\alpha$ we see that $\iota_{n} d\alpha=0$ and hence
$d\alpha=a\, v^\ast\wedge w^\ast$ for some function $a.$ Notice
\begin{align*}
a=d\alpha(v,w) & =-\alpha([v,w])\\
               &=-\frac{1}{\rho}\, g(n,[v,w])=\frac{\theta'}{\rho},
\end{align*}
where the last equality follows from our discussion of $\theta'$ above. Thus
$*d\alpha=\frac{\theta'}{\rho} \, (*v^\ast\wedge w^\ast)
= \frac{\theta'}{\rho}\, n^\ast=\theta'\alpha,$ as claimed.

We now prove that  (2) implies (1).
Again let $v,w$ be a (local) orthonormal frame for
$\xi$ and set $n$ to be the unit normal to $\xi$ for which $v,w,n$ is an
oriented basis for $TM.$ Let $v^\ast,w^\ast,n^\ast$ be the dual basis. We again
see that $\alpha=m\, n^\ast$ for some function $m.$ So 
$*d\alpha=\theta'm\, n^\ast,$ and hence 
$d\alpha = \theta'm\, v^\ast \wedge w^\ast.$ From this we clearly have
$\iota_{n}d\alpha=0.$ Thus $n$ is parallel to the Reeb vector field $R_\alpha$
of $\alpha,$ from which we can conclude that $R_\alpha$ is orthogonal to $\xi.$
\end{proof}

As first observations about weakly compatible metrics we note the following
connections between Riemannian geometry (geodesics) and contact geometry (Reeb
flow lines) and a formula for the mean curvature $H$.

\begin{lemma}\label{lem:covariant-n_n-H}
Let $g$ be a metric that is weakly compatible with a contact form $\alpha$
on the 3--manifold $M.$ Set $n = \frac 1{\rho}R_\alpha$ where 
$\rho = \| R_\alpha\|$ is the norm of the Reeb vector field of $\alpha$.
Then 
\begin{equation}
\label{eqn:nabla_NN_weak}
\nabla_n n= -(\nabla \ln \rho)^\xi,
\end{equation}
where $v^\xi$ denotes the projection of $v$ to $\xi$ along $R_\alpha$.
In particular, $(\nabla \ln \rho)^\xi$ measures the deviation of Reeb flow lines
from tracing out geodesics; and so, if the length of $R_\alpha$ is constant,
then the flow lines of $R_\alpha$ are geodesics. 

In addition we can express the mean curvature of $\ker \alpha$ as
\begin{gather}
\label{eqn:trIIweak}
H = \half \tr \II =-\half \div_{dvol_g} n =-\half n \cdot \left(\ln \frac{\rho}{\theta'}\right).
\end{gather}
\end{lemma}

\begin{proof}	
Since $n$ is a unit vector field, 
we have $\g{\nabla_n n, n} = 0$ and we are left to compute 
$\g{\nabla_n n, v}$ for unit vector fields $v$ tangent to $\xi$. Noting that
$\g{n,\nabla_v n}=0$  and $\rho\,\alpha(v)=\g{n,v}$ (since the two 1--forms in
the formula agree on $\xi$ and $R_\alpha$) we compute 
\begin{align*}
\g{\nabla_n n, v}&=-\g{n,\nabla_n v}=-\g{n, \nabla_n v-\nabla_v n}
= -\g{n, [n,v]} =-\rho\,\alpha([n,v])\\ 
&= \rho\,(d\alpha(n,v)-n\cdot \alpha(v)+v\cdot\alpha(n)) 
= \rho\, (v\cdot\alpha(n))=\rho\, (v\cdot \rho^{-1})\\
&= \rho\, (-\rho^{-2} d\rho(v))=-d(\ln\rho) (v) =\g{-\nabla\ln\rho, v}.
\end{align*}
Thus we see $\nabla_n n= -(\nabla \ln \rho)^\xi.$

In order to compute $H$, we use a (local) moving frame $u, v, n,$ where $u$ and
$v$ are tangent to $\xi.$ We denote by
$\div_\Omega w$ the divergence a vector field $w$ with respect to a volume form
$\Omega$. Recall that $\div_{f\Omega} w= \div_\Omega w+d(\ln f)(w)$ and
$\div_\Omega fw=f\div_\Omega w+df(v).$ Here we have the Riemannian volume form
$\vol_g$  and the contact volume form $\nu = \alpha \wedge d\alpha$ which is
preserved by the Reeb vector field $R_\alpha$. These two volume forms are
related by $dvol_g =\frac{\rho^2}{\theta'} \nu$. Indeed
\[
\alpha\wedge d\alpha(n,u,v)=\alpha(n)d\alpha(u,v)=-\frac 1\rho \alpha([u,v])
= -\frac 1{\rho^2} g(n,[u,v])= \frac{\theta'}{\rho^2}.
\]

Recalling that on
a Riemannian manifold $\div_{dvol_g} n$ can be computed as the trace of the
operator $(u,v)\mapsto \g{\nabla_u n, v}$ we have
\begin{align*}
\tr \II &= \g{n, \nabla_u u + \nabla_v v} 
= - \g{\nabla_u n, u} - \g{\nabla_v n, v} = -\div_{dvol_g} n\\ 
&=  -\div_{dvol_g} (\rho^{-1} R_\alpha) 
= -\rho^{-1}\div_{dvol_g}  R_\alpha - R_\alpha \cdot\rho^{-1} \\
&= -\rho^{-1}\left(\div_\nu  R_\alpha + \theta'/\rho^2 R_\alpha\cdot (\rho^2/\theta')\right) -
R_\alpha\cdot\rho^{-1} \\
&= - \rho^{-2} R_\alpha\cdot \rho + \rho^{-1}/\theta' R_\alpha\cdot\theta' 
= -\rho^{-1}(R_\alpha\cdot\ln \rho - R_\alpha\cdot\ln \theta') 
\end{align*}
So we can conclude that $\tr \II= -n\cdot \ln \frac{\rho}{\theta'}.$
\end{proof}

\no The preceding lemma implies as promised that the two definitions of $\Mg$ in
Equation~\eqref{eq:M_g} agree.

\subsection{Compatibility}
\no We now introduce stronger forms of compatibility.
\begin{definition}
{\em
A contact form $\alpha$ and a Riemannian metric $g$ on a 3--manifold $M$
are \dfn{compatible} if \[
\|\alpha\|=1 \quad \text{and} \quad *d\alpha= \theta' \alpha
\]
for some positive constant $\theta'.$

A contact structure and a Riemannian metric $g$ are compatible if the unit
contact form $\alpha$ of $\xi$ and $g$ are compatible.
(Notice that is is equivalent to saying that the instantaneous rotation of
$\xi$ with respect to $g$ is constant and the Reeb vector field $R_\alpha$ is
unit length and orthogonal to $\xi.$) We may say that $g$ and $\xi$ are
\dfn{strongly compatible} if in addition $\theta'=1$ but we will rarely use
this notion.
}
\end{definition}

Allowing that $\theta'$ be any constant rather than fixing one gives a class of
compatible metrics that is stable under homothety. This is slightly different
from \cite{ChernHamilton85} where $\theta'=2$. Allowing flexibility in $\theta'$
also enables us to absorb various discrepancies stemming from the different
conventions in the definition of the exterior product of differential forms and
their impact on the definition of the exterior derivative. In the current paper,
for any 1--form $\alpha$ and vector fields $u$ and $v$, we use 
$d\alpha(u,v) = u\cdot\alpha(v) - v\cdot\alpha(u) -\alpha([u,v])$. This is twice
the derivative used in \cite{Blair02}.

As a corollary of Lemma~\ref{lem:covariant-n_n-H} we have the
following important geometric properties of compatible metrics.
\begin{corollary}\label{Rderivative}
If $g$ is compatible with $\alpha$ then $R_\alpha$ is a totally geodesic and
divergence free vector field.\qed
\end{corollary}

We now observe that given any metric $g_\xi$ on a contact structure $\xi$ then
there is a canonical way to extend it to all of $M$ so that it is compatible
with $\xi$. This explains in particular the existence of metrics
compatible with a given contact structure (which is well known
\cite{ChernHamilton85}). 
The metric $g_\xi$ induces a volume
form on $\xi.$ Since for any contact form $\alpha_0$ and positive function $f$
we have $d(f\alpha_0)|_\xi=f(d\alpha_0)|_\xi,$ we can find a unique function
$f$ such that $d(f\alpha_0)|_\xi$ agrees with the volume form given on $\xi$
by $g_\xi.$ Setting $\alpha=f\alpha_0$ we can extend $g_\xi$ to $M$ by
demanding that the Reeb vector field $R_\alpha$ is of unit length and
orthogonal to $\xi.$

Note that if $g$ is a metric strongly compatible with $\xi$ then this extension
agrees with the canonical extension of its restriction to $\xi.$ So we have a
``projection'' from the space of metrics $\mathcal{M}(M)$ to
the space of metrics which are strongly compatible with $\xi$. However this
projection can change drastically the geometry of metrics that are compatible
with the  contact structure 
but have $\theta' \neq 1$. We now explain how to remedy this when $M$ is
compact. We denote by $\mathcal{M}_\xi(M)$ the space of metrics compatible with
$\xi$. We want a projection
\[
C: \mathcal{M}(M)\to \mathcal{M}_\xi(M)\subset \mathcal{M}(M).
\]
Given a metric $g$, we denote by $I$ the average instantaneous rotation
speed of $\xi$ with respect to $g$ 
\[
I=\frac{\int_M \theta' \, dvol_g}{\int_M dvol_g}.
\] 
We restrict $g$ to $\xi$ then choose the unique contact form $\alpha$ for which
$\frac 1I d\alpha$ induces the same area form on $\xi$ as $g$ does. Now define 
\[
C(g)= g|_{\xi} + \alpha\otimes\alpha.
\]
If $g$ is compatible with $\xi$ then $I=\theta'.$ One may then use the third
condition of Proposition~\ref{prop:charac_wc} to verify that $C(g)=g.$ 
Moreover, it is clear that $I$ depends smoothly on $g$ and, as the construction
of $C(g)$ from $g$ can be done locally, it is also easy to see that $C(g)$
depends smoothly on $g.$ Thus $C$ is a continuous projection of
$\mathcal{M}(M)$ onto $\mathcal{M}_\xi(M)$.

\section{Comparison of geodesic convexity and symplectic pseudoconvexity}\label{S:convexity_comparison}

In this section we explore various notions of convexity. In 
Subsection~\ref{ss:contactconvexity} we discuss convexity in complex geometry, 
while in Subsection~\ref{ss:riemannianconvexity} we discuss convexity in 
Riemannian geometry. In the last subsection we compare these notions and prove 
Theorems~\ref{thm:compatible}, \ref{thm:weak-compatible} and~\ref{thm:hyp-criterion}.

\subsection{Convexity in almost complex geometry}
\label{ss:contactconvexity}

Let $(W, J)$ be an almost complex manifold.  Assume that $\Omega\subset W$ is an
open, relatively compact, domain in $W.$ Let $\Sigma$ be the hypersurface in $W$
that bounds the region $\Omega$ in $W$ and $f:W\to \R$ be any function for which 
$\Sigma$ is a regular level set, say $\Sigma=f^{-1}(1),$ and $\Omega$ is a
sublevel set of $f.$ The complex tangencies $\mathcal{C}=T\Sigma\cap J(T\Sigma)$
to $\Sigma$ can easily be described as the kernel of the restriction to
$\Sigma$ of the 1--form
$
df\circ J.
$
The form 
\[
L(u,v)=-d(df\circ J)(u,Jv)
\]
is called the \dfn{Levi form} of $\Sigma$.  We say $\Sigma$ is the
\dfn{pseudoconvex}, or \dfn{strictly pseudoconvex}, boundary of $\Omega$ if
$L(v,v)\geq0,$ respectively $L(v,v)>0,$ for all $v\in \mathcal{C}$.

The weak maximum principle for elliptic operators gives a well-known property of
pseudo-convex hypersurfaces: they cannot be ``touched from the inside" by
holomorphic curves. This can be precisely formulated as follows. 
\begin{lemma}\label{th:tangency-lemma}
Suppose $\Sigma = \partial \Omega$ is a strictly pseudo-convex hypersurface in
an almost complex manifold $(W, J)$ as above. Then any 
holomorphic curve $C$ in $\Omega$ 
is transverse to $\Sigma.$ In particular, 
the interior of $C$ is disjoint from $\Sigma.$\qed
\end{lemma}

The almost complex manifolds we will use are symplectizations of contact
manifolds. The symplectization of a contact manifold $(M, \xi)$ is the
submanifold $W$ of $T^*M$ of all covectors vanishing on $\xi$ and inducing the
given co-orientation.  If $\lambda$ denotes the canonical Liouville 1--form on
$T^*M$ then $d\lambda$ is a symplectic form and the fact that $W$ is a
symplectic submanifold of $T^*M$ is equivalent to the contact condition. This
submanifold is a $\R_+$-principal bundle over $M$ and is trivialized by any
choice of section, that is by any contact form $\alpha$. In such a
trivialization, $\R_+\times M$, the symplectic form induced from $T^*M$ becomes
$\omega = d(t\alpha),$ where $t$ is the coordinate on $\R_+$. Let $J$ be a
complex structure on $\xi$ that is compatible with $(d\alpha)|_\xi.$ Extending
$J$ to $TW$ is equivalent to prescribing $J\partial_t$. Gromov's work tells us
that, in order to understand compactness issues for spaces of holomorphic
curves, we should ask $\omega(\cdot, J\cdot)$ to be a Riemannian metric on $W$
(the so called ``tameness condition''). Here this is equivalent to
$\iota_{J\partial_t} d\alpha = 0$ and $\alpha(J\partial_t) > 0$.  This means
$J\partial_t$ should be parallel to the Reeb vector field $R_\alpha$ of
$\alpha$. We will always extend our complex structure from $\xi$ to $TW$ by
setting $J\partial_t=n$ where $n=\frac{1}{\rho}R_\alpha$ and $\rho=\|R_\alpha\|$.

\subsection{Convexity in Riemannian geometry}\label{ss:riemannianconvexity}

Let $S$ be a hypersurface in an Riemannian manifold $(M^n,g)$ that bounds a
region $U$. We say that $U$ is \dfn{geodesically convex at $p$} in $S$ if  any
(local) geodesic in a direction $v\in T_p S$ intersects $U$ only at $p$. The
region $U$ is \dfn{geodesically convex} if it is geodesically convex at every
point $p\in S$. 
\begin{lemma}\label{charicterizeconvex}
Let $f:M^n\to \R$ be a smooth function on a Riemannian manifold $(M^n,g)$ and
let $U$ be a sublevel set of $f$ at a regular value. 
Then $U$ is geodesically convex at $p\in S=\partial U$ if and only if the
Hessian of $f$ is positive definite:
 \[
  \nabla^2 f(v,v) > 0,
 \]
for all non-zero $v\in T_pS.$\qed
\end{lemma}
\no We recall the convexity of the distance function using the functions 
$\textrm{ct}_K$ in Equation~\eqref{cotsdef},  see e.g. \cite{Petersen06}.
\begin{proposition}
\label{prop:riem_convexity}
Let $(M^n,g)$ be a Riemannian manifold and $K$ be a positive constant.
If $g$ has non-positive curvature then geodesic spheres are geodesically 
convex as long as they are embedded. 

If $\sec(g)\leq K$ and
\begin{equation}\label{eqn:convex-radius}
  r < \min\{\inj(g),\frac{\pi}{2\sqrt{K}}\}
\end{equation}
where $\inj(g)$ is the injectivity radius of $(M^n,g),$ then the Hessian
of the distance function 
\[
\mathsf{r}_p:M^n\to \R: q\mapsto d(p,q)
\]
is positive definite on the ball of radius $r$ about $p,$ $B_p(r)$.

More precisely, if $\sec(g)\leq \pm K$, the Hessian of $\mathsf{r}$ satisfies 
\[
\nabla^2\mathsf{r}_p \geq \text{\rm ct}_K(r) g.
\]
\end{proposition}

\subsection{Comparison of convexities}

\no We can now compare the two types of convexity that we considered above.

Let $S$ be a regular level set of a smooth function $f$ on $M$. Let $W$ be the
symplectization of $(M,\xi),$ identified with $\R_+ \times M$ using some
contact form $\alpha$. We set
$\Sigma = \R_+\times S \subset W$ and $\mathcal C
= T\Sigma \cap JT\Sigma$.  Recall that the Levi form of $\Sigma$ is defined on
$\mathcal C$ by $L(u,v) = -d (df \circ J) (u, Jv)$. We want to compare
pseudo-convexity of $\Sigma$, measured by its Levi form $L$, and Riemannian
convexity of $S$, measured by the Hessian $\nabla^2 f.$

\begin{proposition}
\label{prop:levi-form_weak}
Let $g$ be a metric weakly compatible with the contact form $\alpha$ on the
3--manifold $M.$  Then for any
$v\in \mathcal{C}$ we have
\begin{equation*}
\label{eq:Levi-weak} 
L(v,v) =\nabla^2 f(v,v)+\nabla^2 f(Jv,Jv) - \g{\nabla f, \Mg} \|v\|^2  
\end{equation*}
where $\Mg = \nabla_n n + 2H n$, see also Equation~\eqref{eq:M_g}.
\end{proposition}

We use the notation established in Subsection~\ref{ss:contactconvexity}. In particular,
 $J$ will denote the complex structure on the
symplectization $W$ not just on $\xi.$ That is $J$ leaves $\xi$ invariant (and
is the $\pi/2$-rotation of $g_{\|\xi}$) and $J\partial_t=n.$ We also extend the
metric $g$ on $M$ to $W$ by $g+dt\otimes dt.$ We note that the extended $J$ and
$g$ are still compatible.  Using the notation from above, we begin by defining
two bundle maps $A:\xi\to TW$ and $B:\xi\to TW$ by
\[
A(v)=J[Jv, v] - \nabla_v v - \nabla_{Jv} Jv
\]
and
\[
B(v)=J[v, n] + \nabla_{Jv} n + \nabla_n Jv.
\]
One can easily check that $A$ and $B$ are tensorial, meaning that their value at
a point depends on the vector at the point. The main ingredients of  the proof of
Proposition~\ref{prop:levi-form_weak} are the following two technical lemmas.

\begin{lemma}
	\label{lemma:common}
Under the hypotheses of Proposition~\ref{prop:levi-form_weak}, we have
\[
-L(v,v) +\nabla^2 f(v,v)+\nabla^2 f(Jv,Jv) = 
df\left(A(v^\xi) +  B(a Jv^\xi - b v^\xi) - (a^2 + b^2) \nabla_n n \right),
\]
when the vector $v\in \mathcal{C}$ can be written as $v=v^\xi+a\, n+
b\, \partial_t,$ with $v^\xi\in \xi,$ $\R$-invariant and $a$ and $b$ are
constants.
\end{lemma}

\begin{proof}
We first compute
\begin{equation*}
\begin{aligned}
\nabla^2 f (v,v)+\nabla^2 f (Jv,Jv)&= v\cdot (df(v))-(\nabla_v v)\cdot f
+ (Jv)\cdot (df(Jv)) - (\nabla_{Jv}Jv)\cdot f\\
&=v\cdot (df(v)) + (Jv)\cdot (df(Jv)) - df(\nabla_v v +\nabla_{Jv}Jv).
\end{aligned}
\end{equation*}
And, using the formula $d\alpha(u,w) = u\cdot\alpha(w) - w\cdot\alpha(u) -
\alpha([u, w])$ we have
\begin{equation*}
d (df\circ J) (v, Jv) = - v\cdot (df(v)) -(Jv)\cdot (df(Jv)) + df(J [Jv,v]).
\end{equation*}
Adding the two preceding equations, we obtain
\[
-L(v,v) +\nabla^2 f(v,v)+\nabla^2 f(Jv,Jv)= df(J [Jv,v]- \nabla_v v -\nabla_{Jv}Jv).
\]

Decomposing $v$ as $v^\xi+a\, n +b\partial_t$ as in the statement of the lemma
and using $\nabla_{\partial_t} v = 0$ we compute
\begin{gather*}
J[Jv, v] = J[Jv^\xi, v^\xi] + a J[Jv^\xi, n] + b J[n,v^\xi], \\
  \nabla_v v = \nabla_{v^\xi} v^\xi + a(\nabla_{v^\xi} n + \nabla_n v^\xi) +
a^2 \nabla_n n, \\
\intertext{and}
  \nabla_{Jv} Jv = \nabla_{Jv^\xi} Jv^\xi + b(\nabla_{Jv^\xi} n + \nabla_n Jv^\xi)  +
b^2 \nabla_n n.
\end{gather*}
Thus we see that $-L(v,v) +\nabla^2 f(v,v)+\nabla^2 f(Jv,Jv)$ equals
\begin{equation*}
df\left(A(v^\xi) + a B(Jv^\xi) - b B(v^\xi) - (a^2 + b^2) \nabla_n n\right),
\end{equation*}
giving the announced formula.
\end{proof}
\begin{lemma}
	\label{lemma:computeAB_weak}
Under the hypotheses of Proposition~\ref{prop:levi-form_weak}, for any vector
$v$ in $\xi$ we have
\[
A(v) = -\|v\|^2 \left((\tr\II)\, n + \theta' \partial_t\right)
\]
and
\[
B(v) = - \left(\theta' v + (\tr\II) Jv \right) - \g{Jv, \nabla_n n}\, n  
- \g{v, \nabla_n n}\, \partial_t.
\]
\end{lemma}

\begin{proof}
Because $\xi$ is a plane field and both sides of our formulas are tensorial, we
only need to check the announced formulas for a non-zero constant norm vector
field $v$ tangent to $\xi$ and it suffices to check scalar products of our
formula with $v$, $Jv$, $n$ and $\partial_t$.

We first prove that $A(v)$ is normal to $\xi$.  To this end we compute 
\begin{align*}
	\g{A(v), v} &= \g{J[Jv, v], v} - \g{\nabla_v v, v} - \g{\nabla_{Jv}Jv, v} \\
	&= -\g{\nabla_{Jv} v, Jv} + \g{\nabla_{v}Jv, Jv} -\g{\nabla_{Jv}Jv, v}=0.
\end{align*}
One may easily compute, or use the fact that $A$ is $J$--invariant to see, that
$\g{A(v),Jv}=0,$ and thus that $A(v)$ is orthogonal to $\xi.$ The normal
component of $A(v)$ is
\begin{align*}
	\g{A(v), n} &= \g{J[Jv, v], n} - \g{\nabla_v v, n} - \g{\nabla_{Jv}Jv, n} \\
	&= \g{[Jv, v], \partial_t} - \II(v, v) - \II(Jv, Jv) = -\tr \II \|v\|^2.
\end{align*}
Recalling that $\iota_n g=\rho\, \alpha$ and that $\alpha\wedge
d\alpha=\frac{\theta'}{\rho^2} dvol_g$ we compute
\begin{align*}
	\g{A(v), \partial_t} &= \g{J[Jv, v], \partial_t} - \g{\nabla_v v, \partial_t} - 
	\g{\nabla_{Jv}Jv, \partial_t} \\
	&= -\g{[Jv, v], n} =\rho\, d\alpha(Jv,v)=\rho\, \alpha\wedge 
	d\alpha (R_\alpha, Jv, v)\\
	&= -\rho \frac{\theta'}{\rho^2} dvol_g(\rho\, n, v, Jv)= -\theta'\, \|v\|^2.
\end{align*}
The above inner products give the desired formula for $A(v).$

We now compute $B(v)$ by computing its inner product with the chosen
basis. First we have
\begin{align*}
	\g{B(v), v} &= \g{J[v,n],v}+\g{\nabla_{Jv}n,v}+\g{\nabla_n Jv,v}\\
	&= -\g{\nabla_vn, Jv}+\g{\nabla_nv, Jv} +\g{\nabla_{Jv} n, v}-\g{Jv, \nabla_n v,}\\
	&= \g{n, [v,Jv]}.
\end{align*}
As argued above $\g{n,[v,Jv]}=-\theta'\, \|v\|^2.$ Next we have that
\begin{align*}	
\g{B(v), Jv} &= \g{J[v,n],Jv}+\g{\nabla_{Jv}n,Jv}+\g{\nabla_n Jv,Jv}\\
	&=\g{\nabla_vn,v} -\g{\nabla_n v, v} -\g{n, \nabla_{Jv} Jv}\\
	&= -\g{n,\nabla_v v} -\g{n, \nabla_{Jv} Jv} = -(\tr \II)\, \|v\|^2.
\end{align*}
Continuing we compute that 
\begin{align*}
\g{B(v), n}&= \g{J[v,n],n}+\g{\nabla_{Jv}n,n}+\g{\nabla_n Jv,n}\\
  &= \g{[v, n], \partial_t} - \g{Jv, \nabla_n n} = - \g{Jv, \nabla_n n}.\\
\end{align*}
Finally we compute 
\begin{align*}
\g{B(v), \partial_t} &=\g{J[v,n],\partial_t}+\g{\nabla_{Jv}n,\partial_t}
+\g{\nabla_n Jv,\partial_t}\\
  &= -\g{[v, n], n} = -\g{\nabla_v n, n} + \g{\nabla_n v, n}= -\g{v, \nabla_n n},
\end{align*}
These computations yield the stated formula for $B(v).$
\end{proof}

\begin{proof}[Proof of Proposition~\ref{prop:levi-form_weak}]
Given $v\in \mathcal{C}$ since we are only concerned with $v$ at a point 
so we may assume that it is of the form $v=v^\xi+an+b\partial_t$ as in 
Lemma~\ref{lemma:common}. Setting $w = a Jv^\xi - b v^\xi$ 
Lemmas~\ref{lemma:common} and~\ref{lemma:computeAB_weak} yield
\begin{align*}
-L(v,v) +\nabla^2 f(v,v)+\nabla^2 f(Jv,Jv)
&= df\Big(-\tr \II (\|v^\xi\|^2 n + Jw)  
-\theta'(\|v^\xi\|^2 \partial_t + w) \\
&\qquad
- \g{Jw, \nabla_n n}n - \g{w, \nabla_n n}\partial_t
-(a^2 + b^2)\nabla_n n  \Big).\\
\end{align*}

Note that
$w = -bv  + aJv + (a^2 + b^2)\partial_t$ and
$Jw = -av - bJv + (a^2 + b^2)n$.
Thus, keeping in mind that $v$ is in $\mathcal{C}$, we see that $w$ is
in $T\Sigma$ and $Jw$ is in $(a^2 + b^2)n + T\Sigma$. 
Also $\partial_t$ is in $T\Sigma$. So  continuing the above computation we have
\begin{align*}
-L(v,v) &+\nabla^2 f(v,v)+\nabla^2 f(Jv,Jv)\\
&= df\Big(-\tr \II \|v\|^2 n   
- \g{Jw, \nabla_n n}n - \g{w, \nabla_n n}\partial_t
-(a^2 + b^2)\nabla_n n  \Big).\\
\end{align*}
Additionally we notice that 
\begin{align*}
\g{Jw, \nabla_n n}n + \g{w, \nabla_n n}\partial_t &=
\g{v^\xi, \nabla_n n }(v - v^\xi) + \g{Jv^\xi, \nabla_n n }(Jv - Jv^\xi). \\
%
\end{align*}
By the definition of $\mathcal{C}$, $v$ and $Jv$ are in $T\Sigma$, hence in
$\ker df$. So the important part of the right hand-side above is
$-\g{v^\xi, \nabla_n n }v^\xi - \g{Jv^\xi, \nabla_n n }Jv^\xi$ but this is
precisely $ -\|v^\xi\|^2 \nabla_n n$, so we get
the announced formula using the definition $H = \frac 12\tr\II$.
\end{proof}

\no As a direct corollary we obtain the following result. 
\begin{proposition}
\label{prop:levi-form}
Let $g$ be a metric compatible with the contact structure $\xi$ on the
3--manifold $M.$ Then for any
$v\in \mathcal{C}$ we have
\begin{align}
\label{eq:Levi} L(v,v) &=\nabla^2 f(v,v)+\nabla^2 f(Jv,Jv).
\end{align}\qed
\end{proposition}
The later formula was known before when $(W, J)$ is K\"ahler, that is the
Sasakian case, see \cite[Lemma page 646]{GreenWu}.

With these propositions in hand we can prove our main estimates on the tightness
radius of a contact structure. We begin by observing the relation between geodesic 
convexity and pseudo-convexity.

\begin{lemma}
Let $(M,\xi)$ be a contact 3-manifold (not necessarily closed) that is weakly
 compatible with a Riemannian metric $g$ and satisfies $\mg < \infty$. The
 submanifold  $\R_+\times B_p(r)$ of the symplectization of $(M,\xi)$ with the
 complex structure discussed in Section~\ref{ss:contactconvexity} has strictly 
 pseudo-convex boundary if $r<\min\left\{\text{\rm ct}^{-1}_{K}\left(\mg\right), \inj(g)\right\}$.
\end{lemma}
\begin{remark}
Notice that if $\mg=0$, as is the case for compatible metrics, then we see that $\R_+\times \partial B_p(r)$ is pseudo-convex if $\partial B_p(r)$ is geodesically convex. 
\end{remark}

\begin{proof}
We first use Propositions~\ref{prop:levi-form_weak}
and~\ref{prop:riem_convexity} to estimate the Levi form. Denote the radial
function from $p$ by $\mathsf{r}(q)=d(p,q).$ For any complex tangency
$v=v^\xi+a\,n+b\,\partial_t$ of $\R_+\times\partial B_p(r)$ we have
\[
 \begin{split}
L(v,v)& =\nabla^2 \mathsf{r}(v,v) +\nabla^2 \mathsf{r}(Jv,Jv) -
\g{\nabla\mathsf{r}, \Mg}\|v\|^2\\
  & \geq \text{ct}_K(\mathsf{r})(2\,\|v\|^2+a^2+b^2)-\|\Mg\|\|v\|^2 \\
  & \geq (\text{ct}_K(\mathsf{r})- \mg)\| v\|^2 .
 \end{split}
\]
Thus $L(v,v) > 0$ if $r<\text{ct}_K^{-1}(\mg)$, 
\end{proof}

\begin{proof}[Proof of Theorem~\ref{thm:compatible} and Theorem~\ref{thm:weak-compatible}]
Let $B_p(r)$ be the geodesic ball of radius $r$ about $p$. From the last lemma  we know that 
for any 
\[
r <\min\left\{\text{\rm ct}^{-1}_{K}\left(\mg\right), \inj(g)\right\},
\]
$\R_+\times\partial B_p(r)$ is a
strictly pseudo-convex submanifold of $\R_+\times M$. 
Assume now for contradiction that $B_p(r)$ contains an overtwisted disk.
Arguing as in \cite{Hofer93}, one can start a family of holomorphic disks 
near an elliptic singularity in the overtwisted disk. Since these disk cannot touch 
$\R_+\times\partial B_p(r)$ thanks to pseudo-convexity, the Gromov--Hofer
compactness theorem extends to this setting and we get the existence of a 
closed Reeb orbit $\gamma$ in 
$B_p(r)$ and a $J$-holomorphic cylinder $C_\gamma=\R_+\times\gamma$ in
$\R_+\times B_p(r)$. But this is a contradiction because $C_\gamma$ has
to touch $\R_+\times\partial B_p(r_0)$ from the inside for some $r_0\leq r$.
\end{proof}
\begin{proof}[Proof of Theorem~\ref{thm:hyp-criterion}]
Note that pull-backs of $\xi$ and the metric to any covering space are (weakly)
compatible and the sectional curvature is non-positive. It is well known, by
Hadamard Theorem \cite{Chavel06}, that the
universal cover of a three manifold with nonpositive curvature is $\R^{3}$ 
and the space is exhausted by geodesic balls.
The assumption $\mg\leq \sqrt{K}$
implies
\[
 \text{ct}_K(r)- \mg > 0,\qquad \text{ for all $r\geq 0$},
\]
since $\text{ct}_K(r) > \sqrt{K}$.
Theorem~\ref{thm:weak-compatible} says a ball of any radius is tight, but since
any potential overtwisted disk will have to be contained in a ball of some
radius it cannot exist.
\end{proof}

\section{Tightness and quasi-geodesics}\label{S:quasi-geod}
Here we establish a geometric ``universal tightness'' criterion for contact
structures even when the metric and contact structure are not compatible. 

\begin{proof}[Proof of Theorem~\ref{thm:vwc}]

Suppose by contradiction $(M,\xi)$ is overtwisted.  According to \cite{Hofer93}
there
will be a closed contractible orbit in the flow of the Reeb field $R,$ and hence
in the flow of $N=R/\|R\|.$ This orbit will lift to a closed orbit $\gamma$ in
the universal cover of $M.$ Of course our metric, contact structure and Reeb
field also lift to the universal cover and satisfy the same hypotheses (since we
will work exclusively in the cover from now on, we will use the same notation
for objects in the cover).

Choose any point $p$ in the cover. There will be some $r$ such that the embedded
geodesic ball $B_p(r)$ of radius $r$ about $p$ will contain $\gamma$ and
$\partial B_p(r)$ will have a tangency with $\gamma.$ Let ${\mathsf
r}_p(x)=d(p,x)$ be the radial function measuring the distance from $p.$ The
estimate from Proposition~\ref{prop:riem_convexity} says
\[
\sqrt{K} < \nabla^2 {\mathsf r}_p (\dot\gamma, \dot\gamma).
\]
Since we are assuming $\gamma$ has tangency to $\partial B_p(r)$ we can
parameterize $\gamma$ such that the tangency occurs at $\gamma(0).$ As $\gamma$
lies inside $B_p(r)$ we see that 
\[
0\geq   \frac{\partial^2}{\partial t^2}({\mathsf r}_p\circ \gamma)|_{t=0} 
= \frac{\partial}{\partial t} g(\nabla {\mathsf r}_p, \dot \gamma)
= \nabla^2{\mathsf r}_p(\dot\gamma,\dot\gamma)+ d{\mathsf r}_p 
(\nabla_{\dot\gamma}\dot\gamma).
\]
Thus we can compute
\[
\sqrt{K} < \nabla^2 {\mathsf r}_p (\dot\gamma, \dot\gamma)\leq - 
g(\nabla{\mathsf r}_p,\nabla_{\dot \gamma}\dot\gamma)\leq
 \|\nabla_{\dot\gamma}\dot\gamma\|=\|\nabla_N N\| \leq \sqrt K,
\]
where the third inequality comes from the fact that $\nabla {\mathsf r}_p$ has
unit length so $ - g(\nabla{\mathsf r}_p,\nabla_{\dot \gamma}\dot\gamma)$ is one
component of $\nabla_{\dot\gamma}\dot\gamma$ in some orthonormal basis. Thus we
arrive to the absurd consequence that $\sqrt K < \sqrt K$ and an overtwisted
disk could not exist.
\end{proof}

\begin{remark}
{\rm 
 The above proof works also in the higher dimensional setting and shows that
the condition of Theorem~\ref{thm:vwc} implies the contact manifold is plastikstufe free 
(see \cite{Niederkrueger06} for a definition of  plastikstufe, and \cite{EKM-higher} for further discusion.)
}
\end{remark}

\section{Characteristic foliations and contact topology}
\label{S:xiS}

In this section we recall the definition of the characteristic foliation of a
surface in a contact 3--manifold and use it to study the tightness radius of a
compatible metric. 
\subsection{Characteristic foliations}
Recall that an oriented singular (this adjective will be implicit in the
following) foliation on an oriented surface $S$ is an equivalence class
of 1--forms on $S$ where $\alpha \sim \beta$ if there is a positive function $f$ such
that $\alpha = f\beta$. Let $\alpha$ be a representative for a singular
foliation $\F$. A singularity of $\F$ is a point where $\alpha$ vanishes. The
singularity $p$ is said to have non-zero divergence if
$(d\alpha)_p$ is an area form on $T_pS.$ If $\omega$ is
an area form on $S$ (compatible with the chosen
orientation) then to 
each singular point $p$ we attach the sign of the unique real number $\mu$ such that 
$(d\alpha)_p = \mu \omega_p$. One can easily check that singular points and their signs
do not depend on the choice of $\alpha$ in its equivalence class or on $\omega$ if
we keep the same orientation.

Let $S$ be an oriented surface in a contact manifold $(M, \xi)$ with $\xi = \ker
\alpha,$ co-oriented by $\alpha$. The characteristic foliation $\xi S$ of $S$ is
the equivalence class of the restriction of $\alpha$ to $S$. The contact
condition ensures that all singularities of characteristic foliations have
non-zero divergence. Singularities of $\xi S$
correspond to points where $S$ is tangent to $\xi$ and they are positive or
negative according as the orientation of $\xi$ and $S$ match or do not match.
We also notice that $\alpha$ provides a co-orientation, and hence if $S$ is 
oriented by an area form $\omega$ the orientation of the line field $\xi S$ is given by 
the vector field $X$ which satisfies $\iota_X \omega=\alpha|_S$.

One may dually think of the characteristic foliation on $S$ as coming from the
singular line field on $S$ given by $T_pS\cap \xi_p$ for each $p\in S.$

\subsection{Convex spheres and characteristic foliations on spheres}

An oriented foliation $\F$ on a sphere $S$ is \dfn{simple} if it has exactly one
singularity of each sign (the positive one will be called the north pole and the
negative one the south pole).  It is \dfn{almost horizontal} if, in addition,
all its closed leaves are oriented as the boundary of the disk containing the
north pole (in other words from ``west'' to ``east'').  These are (slight
variations on) definitions due to Eliashberg \cite{Eliashberg_vrille}.

If $\xi$ is a contact structure on $M$ and $S$ is a sphere in $M$ then 
we will say that $\xi$ is simple or almost horizontal along $S$ if $\F = \xi S$
has this property. The relevance of these definitions to compatible metrics
comes from the following lemma.

\begin{lemma}
\label{lem:metricsimple}
Let $g$ be a Riemannian metric compatible with the contact manifold $(M,\xi).$
Let $\alpha$ be the contact form implicated in the definition of compatibility
between $g$ and $\xi$ and $R_\alpha$ its Reeb vector field.  Let $r<\inj(g)$ and
$S$ be the sphere of radius $r$ around some point $p_0.$ The contact structure
$\xi$ is simple along $S$ with poles $\exp_{p_0}(\pm rR)$.
\end{lemma}

\begin{proof}
Let $S$ be a sphere of radius $r < \inj(g)$ around $p_0.$ Let $\gamma$ be a geodesic
starting at $p_0$ and denote $\gamma(r)\in S$ by $p$. Suppose $p$ is a singularity of
$\xi S,$ that is, $\xi_p$ is tangent to $S$ at $p.$ By Gauss' Lemma, $\xi_p$ is
orthogonal to $\gamma$ at $p.$ As $R_\alpha$ is also orthogonal to $\xi_p$ we
must have $\gamma'(r)=\pm R_\alpha.$ As the flow of $R_\alpha$ give geodesics we
see that $\gamma'(t)$ is equal to $\pm R_\alpha$ for all $t$ and hence
$\gamma'(0)=\pm R_\alpha.$ Thus $p$ is $\exp_{x_0}(\pm rR_\alpha)$.
\end{proof}

We now have the following important definition in contact geometry due to Giroux
\cite{Giroux91}. 
We say that an hypersurface $S$ in a contact manifold $(M,\xi)$ is
\dfn{$\xi$--convex} if there is a vector field transverse to $S$ whose flow preserves
$\xi$. (Normally such a hypersurface is just called convex, but to distinguish
this notion from other types of convexity in this paper we use the term
$\xi$--convex.) In our situation, we can recognize $\xi$--convex spheres using
the following very special case of a criterion by Giroux.

\begin{lemma}[Giroux 1991, \cite{Giroux91}]
\label{lem:crit_convexity}
If a contact structure $\xi$ is simple along a sphere $S,$ then $S$ is
$\xi$--convex if and only if $\xi S$ has no degenerate closed leaf.\qed
\end{lemma}

The following proposition explores how a contact structure on a ball can be
overtwisted by explaining some relations between simple, almost horizontal and
tight. The first point is obvious while the second one was observed by Giroux.
The third one will be crucial for the sphere theorem. The second point will not
be directly used in this paper but it could prove useful to get tightness radius
estimates in later work and its proof is also needed for the third point.

\begin{proposition}
	\label{prop:horiztight}
Let $B$ be a ball in a 3--dimensional contact manifold $(M, \xi)$ which is the
disjoint union of a point $p$ and a family of spheres $S_t$, $t \in (0, 1].$ 
\begin{enumerate}
\item 
If all foliations $\xi S_t$ are simple and $\xi|_{B}$ is tight then $\xi$
is almost horizontal along all the $S_t.$
\item 
If all foliations $\xi S_t$ are almost horizontal, then $\xi|_{B}$ is tight.
\item 
If all foliations $\xi S_t$ are simple, and $\xi|_{B}$ is overtwisted then
there is some radius $t_1$ such that all foliations $\xi S_t$ for $t \geq t_1$
have closed leaves.
\end{enumerate}
\end{proposition}

\begin{proof}
To establish the first point we notice that  tightness rules out the existence
of any closed leaf in any of the characteristic foliations $\xi S_t,$  since
such a leaf would bound an overtwisted disk. So they are all (trivially) almost
horizontal.

We now prove a special case of the second point of the proposition: if there
is no closed leaf in any $\xi S_t$ then $\xi$ is tight.
Indeed, in the absence of closed leaves, Lemma~\ref{lem:crit_convexity} ensures
that all the spheres $S_t$ are $\xi$--convex.
The dividing curve $\Gamma_t$ on each $S_t$ can be taken to be any simple closed
curve separating the north and south pole and transverse to $\xi S_t$, so we
can assume the $\Gamma_t$ depend smoothly on $t.$ We can then construct a
contact vector field that is transverse to all the $S_t$ (this can be done by
using the fact that the space of contact vector fields supported in a
neighborhood of a convex surface $S$ is contractible and the fact that such a
vector field is transverse to all surfaces close to $S$).  We can use this
vector field to embed $(B, \xi)$ in the standard contact structure on $\R^3$.
Alternatively, once we know that all spheres $S_t$ are $\xi$--convex, we can
directly apply Giroux's reformulation of Bennequin's theorem, 
\cite[Theorem~2.19]{Giroux00} to get tightness.

We now assume that closed leaves do appear.
Darboux's theorem, coupled with Bennequin's result that the standard contact
structure on $\R^3$ is tight, ensures that there is some $\varepsilon > 0$ such 
that spheres $\xi S_t$ for $t \leq \varepsilon$ have no closed leaf. 
So there is a minimal $t_1 > 0$ such that $\xi S_{t_1}$ has a closed leaf.
In order to finish the proof of the second point of the proposition, we have to
prove that such a leaf has to go to the west ({\em i.e.}~it is oriented as the
boundary of the disk containing the south pole), contradicting the almost
horizontal assumption.
In order to prove the third point, we have to prove that all spheres $S_t$ for
$t \geq t_1$ exhibit closed leaves.
This will be done using essentially the analysis in Giroux's Birth-Death Lemma, 
\cite[Lemma~2.12]{Giroux00}.
In order to facilitate this analysis, we observe the following result. 
\begin{lemma}\label{technicallemma}
Let $A$ be an annulus. 
There is an embedding of $A\times[\varepsilon, 1]$ into $B$ such that 
\begin{itemize}
\item[$(i)$]
each $A_t := A \times \{t\}$ maps into the corresponding $S_t$,
\item[$(ii)$]
the complement of $A_t$ in $S_t$ is made of two disks $D_t^n$ and $D_t^s$
(around the poles) whose characteristic foliations is topologically conjugated
to a radial foliation and goes transversely out of $D_t^n$ and into $D_t^s$
(recall these foliations can be oriented),
\item[$(iii)$]
the vector field $\partial_t$ on $A \times [\varepsilon, 1]$ maps to a
Legendrian vector field,
\end{itemize}
Moreover, by a $C^\infty$ small perturbation of 
$A\times[\epsilon,1]$ that is fixed near any
finite number of $A_t$ for which the foliation $\xi A_t$ is non-degenerate, we may assume in
addition that
\begin{itemize}
\item[$(iv)$] for $t$ outside a finite set $\varepsilon <t_1 < \dots < t_N < 1$, all
	closed leaves of $\xi S_t$ are non-degenerate and hence stable, and
\item[$(v)$] 	each $\xi S_{t_i}$ has exactly one degenerate closed leaf and it 
	indicates the birth or death of a pair of non-degenerate closed leaves.
\end{itemize}
\end{lemma}
We will establish this technical lemma (which belongs purely to the realm of
dynamical systems) after we finish the proof of our proposition.

The pulled back contact structure on $A \times [\varepsilon, 1]$ will also be
denoted by $\xi$.
The boundary component of each the annuli $A_t$ (or their essential sub-annuli
below) next to the disk $D_t^n$ (respectively $D_t^s$) will be called
the northern (respectively southern) border. Similarly, a closed leaf in any 
$\xi A_t$ will be called northernmost if it bounds an open disk containing
$D_t^n$ and no closed leaf.

The following lemma is a variation on Giroux Birth-Death Lemma.
Here we see the contact condition with specific orientation plays a crucial
role by determining the type (birth or death) of bifurcations as $t$ increases.
\begin{lemma}\label{auxlem}
Let $A' \subset A$ be a subannulus of $A$ containing the northern border
and pick $\bar t$ in $(\varepsilon, 1)$. Suppose $\xi A'_{\bar t}$ has a
closed leaf. If there is some positive $\delta$ such that $\xi A'_t$ has no
closed leaf when $\bar t < t < \bar t + \delta$ (respectively 
$\bar t - \delta < t < \bar t$) then the northernmost closed leaf of
$\xi A'_{\bar t}$ is oriented from west to east (respectively east to west).
\end{lemma}

\begin{proof}
We denote by $\gamma$ the northernmost closed leaf of $\xi A'_{\bar t}$.
Suppose it goes to the west (the other case being symmetric) and let $A''$ be
the sub-annulus of $A'$ bounded by $\gamma$ and the northern border of $A$.

Because $\xi$ is a positive contact structure and $\partial_t$ is Legendrian,
the characteristic foliations on the annuli $A_t$ are rotating strictly
clockwise at each point of $A$. We know the flow is pointing in along the
northern border of $A''$. Moreover it is tangent to $\gamma$ at $t = \bar t$.
Thus our assumption on the direction of the flow along $\gamma$ implies there
is some positive $\delta$ such that the characteristic foliation of $A''_t$
will point inward along $\gamma$ for $\bar t< t< \bar t + \delta$. So, for
these $t$, $\xi A''_t$ points inward along both boundary components. Since 
$\xi A''_t$ is non singular, the Poincar\'e-Bendixson theorem then guaranties
the existence of an interior closed leaf.
\end{proof}

The second point of the proposition above clearly follows from the above lemma
with $A' = A$ and $\overline{t}= t_1$.
The lemma also proves that births next to the northern boundary are always
births of leaves heading west (here $A'$ is a annulus containing the new closed
leaf but no other one at the birth time).  The same reasoning with $t$ going
backward shows that deaths next to the northern boundary are always deaths of
leaves heading east.

Let us now complete the proof of the the third point. Assume for contradiction
that closed leaves disappear completely at time some time $t'$. Considering the
region $S\times [\varepsilon, t']$ we may use Lemma~\ref{technicallemma} 
to find the annulus $A$ and $A\times[\varepsilon,t']$ satisfying the $(i)$--$(iii)$. Notice 
that $A_\varepsilon$ and $A_{t'}$ have non-degenerate characteristic foliations
as there are no singularities or closed leaves. Thus we may $C^\infty$-perturb 
$A\times[\varepsilon,t']$ relative to $A_\varepsilon\cup A_{t'}$ to satisfy $(iv)$--$(v)$ too.
Let $\varepsilon<t_1<\ldots, t_k<t'$ be the times from item $(iv)$ in the lemma. 
We know by the 
preceding paragraph that the (unique and northernmost) closed leaf at $t_1$ goes to the
west whereas the one at $t_k$ goes to the east. In addition there are finitely
many births and deaths of closed leaves, they occur at times $t_i$. We now
prove by induction on these events that the northernmost closed leaf always
goes to the west, contradicting what we just mentioned at $t_k$. We already
know this at birth time $t_1$. Any birth or death not affecting the
northernmost orbit is irrelevant. Births next to the northern border spawn
leaves going to the west. Finally there cannot be any deaths next to the
northern boundary because they would involve a northernmost orbit heading east
but they are prohibited by induction hypothesis.  
\end{proof}

We now turn to the proof of our technical Lemma~\ref{technicallemma}.
This proof was kindly provided by Fr\'ed\'eric Le Roux. 
\begin{proof}[Proof of Lemma~\ref{technicallemma}]
First observe that the poles are both index one singularities. This can be
checked directly for very small radii and then follows by continuity since we
know that all foliations $\xi S_t$ have exactly two isolated singularities.
In addition these singularities have non-zero divergence because they come from
a contact structure. The first important point concerns the picture near these
singularities.  

\vspace{.5cm}
\noindent
{\bf Claim: }Let $p$ be an isolated singularity of a vector field $X$ on a
surface. If $X$ has non-zero divergence at $p$ and index one then there is a
smooth disk $D$ around $p$ such that $X$ is transverse to the boundary of $D$
and defines a foliation which is topologically conjuguated to a radial
foliation in $D$. In addition, if $X$ belongs to a smooth one-parameter family
of such vector fields then $D$ can be chosen to depend smoothly on the
parameter as well.

\begin{proof}[Proof of the claim]
Let $D_0$ be a disk around $p$ where the divergence of $X$ (for some
fixed auxilliary area form) does not vanish. Shrinking $D_0$ if necessary, we
can assume that $p$ is the only singular point in $D_0$. Poincar\'e--Bendixson's
theorem then tells us there are only two kinds of orbits of $X$ which are
entirely contained inside $D_0$ (besides $p$ itself): closed orbits and orbits
whose $\alpha$ and $\omega$--limits are both $p$. In both cases there would
exist a subdisk $D$ of $D_0$ such that the flux of $X$ through $\partial D$
vanishes. But Stokes' theorem shows that this cannot happen under our non-zero
divergence assumption. So there is no nontrivial orbit contained in $D_0$.

Now let $\gamma$ be a smooth simple closed curve going once around $p$ in $D_0$
and which minimizes the number of tangencies with $X$ among such curves. We will
prove that $\gamma$ cannot have any tangency with $X$ hence the disk $D$ bounded
by $\gamma$ is the disk we sought. A priori, tangencies with $X$ come in two
flavors. The orbit of $X$ through the tangency point can be locally inside or
outside the interior of $D$. But the minimizing property of $\gamma$ forbids
interior tangencies.  Indeed, if $q$ is such a tangency point then the orbit of
$X$ through $q$ has to leave $D$ at some point $q'$  because it cannot be
contained in $D$. Taking a long flow box along the piece of orbit between $q$
and $q'$, we see that we can modify $\gamma$ to reduce the number of tangency
points, a contradiction.  Hence $\gamma$ has only exterior tangencies with $X$.
But the index of $X$ at $p$ is one plus half the number of interior tangencies
minus half the number of exterior ones so, using the index assumption, we get
that $X$ is transverse to $\gamma$. To get the topological picture inside $D$,
it suffices to remark that the Poincar\'e-Bendixson theorem now guaranties that
all trajectories in $D$ go from the boundary to $p$.

We now want this to work for one-parameter families. Any disk $D$ satisfying 
the conclusions of the first part with respect to some vector field $X_t$ in
our family also satisfies them for $X_{t'}$ when $t'$ is close to $t$ because
 transversality to the boundary of $D$ is an
open condition. So we only need to prove that, for a fixed $X$, the space of
suitable disks is path connected. Let $D$ and $D'$ be two such disks. Any
trajectory through the boundary of $D$ hits the boundary of $D'$ exactly once.
So we can build an isotopy pushing $D$ to $D'$ following these trajectories.
\end{proof}

We now return to our family of simple spheres $S_t$. We denote by $D_t^n$ and
$D_t^s$ two families of disks given by the above claim around the north and
south poles respectively.
Notice that each $S_t\setminus (D_t^n\cup D_t^s)$ is an annulus transverse to
$\xi$.  Thus if $A$ is an annulus then we can find an embedding of
$A\times[\varepsilon, 1]$ such that (i) through (iii) from the lemma hold. 

Since the image of this embedding is compact and the characteristic foliation
on $A_t$ for all $t$ is non-singular we know that the maximal angle between
$\xi_p$ and $T_pA_t$ is bounded below. Thus we can perturb the $A_t$ so as to
make the characteristic foliations generic without introducing any
singularities. The only thing we want to check is that we can make this
perturbation without getting rid of all closed leaves. So we proceed in two
steps. First we perturb the family close to any chosen closed leaf to make that
leaf non-degenerate. Then we perturb the whole one-parameter family of
foliations on $A$ like in \cite{Sotomayor} using a small enough perturbation so
that the stable leaf created during the first perturbation does not disappear.
\end{proof}

Note that the hypotheses of the second point in the proposition are trivially 
satisfied if no characteristic foliation of a sphere $S_t$ has a closed leaf. This 
leads to the remarkable fact that if $\xi$ is overtwisted on $B$ then this can be 
seen clearly on all spheres $S_t$ for sufficiently large $t$.

\begin{proof}[Proof of Theorem~\ref{thm:onspheres}]
Fixing $p\in M$ consider the geodesic spheres $S_p(r)$ of radius $r$ about $p$
and the geodesic balls $B_p(r)$ that they bound. We can use
Lemma~\ref{lem:metricsimple} to conclude that all the spheres $S_p(r), r\leq
\inj_p(g),$ have simple characteristic foliation. Recall that we are assuming
that $\tau_p<\inj_p(g)$ {\em i.e.} $B_p(r),r< \tau_p,$ is tight and $B_p(r)$,
$r>\tau_p,$ is overtwisted.  Let 
\[
r'=\inf\{r\, |\, \text{such that $S_p(r)$ has
a closed leaf in its characteristic foliation}.\}
\]
 Notice that $S_p(r')$ does
have a closed leaf since simple foliations on spheres without closed leaves form
an open set. By Proposition~\ref{prop:horiztight} the contact structure
restricted to $B_p(r),r<r',$ is tight. Thus $r'=\tau_p$ and we see that
$S_p(\tau_p)$ has a closed leaf  in its characteristic foliation, which of
course bounds an overtwisted disk. We then get overtwisted disks on spheres of
higher radii using the third point of Proposition~\ref{prop:horiztight}.
\end{proof}

\section{The contact sphere theorem}\label{S:CST}

In this section we prove the contact sphere theorem and discuss possible generalizations of it. 

\subsection{Proof of the contact sphere theorem}
The following proposition is a variation on a similar result used in 
the proof of the topological sphere theorem \cite{Klingenberg61}, it does not involve 
any contact geometry.

\begin{proposition}\label{prop:pinching}
Let $M$ be complete simply connected Riemannian manifold whose sectional
curvature satisfies $\frac{4}{9} < K \leq 1$. If $p$ and $q$ in $M$ are at
maximal distance, that is $d(p, q) = \diam(M)$, then there are radii $r_p$ and
$r_q$ such that 
\begin{itemize}
\item 
both closed balls $\overline{B}(p, r_p)$ and $\overline{B}(p, r_p)$ are embedded, {\em i.e.}
$r_p, r_q < \inj(g)$,
\item 
the ball $\overline{B}(q, r_q)$ is convex, {\em i.e.} $r_q < \conv(g)$,
\item 
$M = \overline{B}(p, r_p) \cup \overline{B} (q, r_q)$, and
\item
the boundary of each ball, $\overline{B}(p, r_p)$ and $\overline{B}(p, r_p)$, is contained in the interior of the other ball.
\end{itemize}
\end{proposition}

\no We remark that the Bonnet Theorem guaranties that $M$ is compact in the above
proposition so points like $p$ and $q$ do exist.

\begin{proof}
Because $M$ is compact, there is some $\delta$ such that
$\frac{4}{9} < \delta < K \leq 1$. Let $\varepsilon_p$ and $\varepsilon_q$ be
positive numbers  to be fixed later.
We set $r_p = \pi(1 - \varepsilon_p)$ and $r_q = \frac{\pi}{2}(1 - \varepsilon_q)$.
Klingenberg's injectivity radius estimate, see \cite[Theorem~5.10]{CheegerEbin}, is
$\inj(g) \geq \pi$ so that both balls are embedded. The convexity radius is
then at least $\pi/2$, see \cite[Theorem~7.10]{Chavel06}, so our ball around $q$
is convex. To prove the remaining two properties, it suffices to prove that,
for any $x$ in $M$, $d(x, q) \geq r_q$ implies that $d(x, p) < r_p$.

Let $\gamma_1$ be a geodesic segment between $q$ and $x$ so that 
$L(\gamma_1) = d(q, x) > r_q$. According to Berger's lemma
\cite[Lemma~6.2]{CheegerEbin}, there is a geodesic segment $\gamma_2$ from $q$
to $p$ such that the angle $\alpha$ between $\gamma_1$ and $\gamma_2$ at $q$ is
not greater than $\pi/2$. We want to get an upper bound on $d(p, x)$. According
to Toponogov theorem, we only need to get it for the corresponding hinge in
the round 2--sphere in Euclidean space with curvature $\delta$, hence radius 
$R = \frac{1}{\sqrt\delta}$.

Rather than considering lengths, it is more convenient to work with the angle
under which a segment is seen from the center of the sphere. The actual length
is $R$ times this angle. Let $a$, $b$ and $d$ be angles such that the triangle
we are looking at has sides whose length
are $Rb$ coming from $\gamma_1$, $Rd$ coming from $\gamma_2$ and $Ra$ which
is not less than $d(p, x)$. All these angle are between $0$ and $\pi$.
The spherical law of cosines tells us
\[
\cos(a) = \cos(b)\cos(d) + \sin(b)\sin(d)\cos(\alpha)
\]
The second term is non-negative because $\alpha \leq \pi/2$. Since $x$ is not
in the ball around $q$, we have $Rb \geq \pi(1-\varepsilon_q)/2$, so 
$\cos(b) \leq \cos(\pi\frac{1-\varepsilon_q}{2R})$ and 
and we get, from the law of cosines above, 
$a \leq \pi - \pi\frac{1-\varepsilon_q}{2R}$.
The corresponding distance is $\pi\left(R - \frac{1-\varepsilon_q}{2}\right)$ 
which is strictly less than $\pi(1 - \varepsilon_p) = r_p$ whenever 
$2\varepsilon_p + \varepsilon_q < 3 - \frac{2}{\sqrt\delta}$. The lower bound
hypothesis on $\delta$ allows to choose positive $\varepsilon_p$ and
$\varepsilon_q$ satisfying this inequality.
\end{proof}

\begin{proof}[Proof of Theorem~\ref{CST-thm}]
We now gather the different ingredients of the contact sphere theorem,
highlighting how Riemannian geometry, topological methods in contact geometry
and pseudo--holomorphic curves arguments interact in this proof.

Since both the hypotheses and the conclusion of the theorem are scale
invariant, we can assume that the curvature is bounded above by $1$ so that
$\frac{4}{9} < K \leq 1$. Moreover, we can assume that $M$ is simply connected
as pulling  the contact structure and metric back to the universal cover of $M$
does not affect the curvature pinching. 

Deep classical Riemannian geometry gives, through
Proposition~\ref{prop:pinching}, that there are two geometrics balls
$B_\text{cvx}$ and $B_\text{big}$ whose interior covers our manifold $M$ and
such that $B_\text{cvx}$ is weakly convex.

We assume for contradiction that $\xi$ is an overtwisted contact structure.  A
priori, it could be that all overtwisted disks intersect both $B_\text{cvx}$
and $B_\text{big}$. However there is no loss of generality in assuming that
there is one, which we denote by $D$, which misses the center $q$ of
$B_\text{cvx}$.

Our comparison of Riemannian and almost-complex convexity combines with
pseudo-holomorphic curves arguments of Gromov and Hofer to tell us, through 
Theorem~\ref{thm:compatible}, that $B_\text{cvx}$ is a tight ball.

We can now use either the classification of tight contact structures on balls
by Eliashberg \cite{Eliashberg92a} or, more elementarily, our description in
Proposition~\ref{prop:horiztight} and its proof. Either way, we can construct a
contact vector field transverse to the concentric spheres that make up
$B_\text{cvx} \setminus \{q\}.$ This contact vector fields generates a contact
isotopy that will push any subset of $B_\text{cvx}$ that misses $q$ into an
arbitrarily small neighborhood of $\partial B_\text{cvx}.$ In particular we can
push $D$ into $B_\text{big}$.

Based on the special geometry of compatible metrics and a topological argument
using Giroux's study of bifurcations for characteristic foliations,
Proposition~\ref{prop:horiztight} then tells us that there is an overtwisted
disk on $\partial B_\text{big}.$ However, as we know $\partial B_\text{big}\subset
B_\text{cvx},$ this contradicts the tightness of $B_\text{cvx}.$ Hence we
see that $\xi$ must be tight.  

Although this is does not completely follow from the previous discussion, the
ambiant manifold is the 3-sphere as is guarantied by the classical sphere
theorem. Now we get that $\xi$ is isomorphic to the standard contact structure
because all tight contact structures on the sphere are standard. This later fact is due
to Eliashberg \cite{Eliashberg92a} and uses purely topological method in
contact geometry (see also \cite[Remark~2.18]{Giroux00} for Giroux's alternative
proof).  
\end{proof}

\begin{remark}\label{rem-GG}
{\em
As an immediate corollary of the contact sphere theorem, we obtain that any
contact structure compatible with the round metric on $\S^3$ is isotopic to the
standard one $\xi_0$.
While this result is not obvious, it already follows from older results.
Indeed, suppose that $\xi$ is compatible with a round metric and denote by $R$
the Reeb vector field involved. Corollary~\ref{Rderivative} guaranties that $R$
is geodesic (its orbits are great circles parametrized by arc length) and
divergence free.
Although there is an infinite dimensional space of geodesic vector fields on
$\S^3$, Gluck and Gu proved in \cite{GluckGu01} that they become completely
rigid if we assume in addition that they are divergence free.
So $\xi$ is actually conjugated to $\xi_0$ by an isometry.

This rigidity is of a completely different nature than the statement of the
sphere theorem where the isotopy with the standard structure is unrelated to
any rigid structure on the sphere. 
}
\end{remark}

\subsection{Extensions of the contact sphere theorem.}\label{sec:extensions}
We now discuss possible extensions and strengthening of the contact sphere
theorem. We first remark that the hypothesis on the metric is open in the space
of compatible metrics. This space can be seen as the product of the space of metrics
on $\xi$ and the space of positive constants $\theta'$. We now describe one
way to alter the statement of Theorem~\ref{CST-thm} so that it applies to an
open set in the space of all metrics. It uses the projection $C$ defined at the
end of Section~\ref{S:compatible_metrics}.

Given any metric $g$ on $M$ we define the \dfn{$\xi$-adapted sectional
curvature} for $g$ to be the sectional curvature of $C(g),$ and denote this
curvature as $\sec_\xi(g).$ We have the following immediate corollary of
Theorem~\ref{CST-thm} and the fact that the projection $C$ depends smoothly on
$g.$

\begin{corollary}
Let $(M,\xi)$ be a closed contact manifold and $g$ any complete Riemannian
metric. If there is a constant $K_{max}>0$ such that the $\xi$-adapted
sectional curvatures of $g$ satisfy 
\[
0<\frac 49 K_{max}<\sec_\xi(g)\leq K_{max},
\]
then the universal cover of $M$ is diffeomorphic to the 3--sphere by a
diffeomorphism taking the lift of $\xi$ to the standard contact structure on
the 3--sphere. 

Moreover, the set of metrics satisfying the the $\xi$-adapted sectional
curvature pinching is open in the set of all metrics. 
\end{corollary}

It would be interesting to see how the sectional curvature and the
$\xi$-adapted sectional curvature of $g$ compare. For instance,  there should
be a version of the sphere theorem involving weakly compatible metrics, where
the pinching constant is weakened by a factor depending on the derivatives of
$\rho$ and $\theta'$.

Another natural question is if the pinching constant $4/9$ in Theorem~\ref{CST-thm} 
is optimal. Recall that  the classical pinching constant $1/4$ of
the sphere theorem is optimal in even dimensions. Moreover, in all dimensions,
the classical proof of  the sphere theorem shows that any simply connected
complete Riemannian manifold whose sectional curvature is $1/4$-pinched can be
covered by two balls.  In trying to extend this idea to the contact geometric
setting, notice that the techniques of Eliashberg mentioned above allow one to
prove that a contact 3--manifold which can be covered by two standard balls is
tight, hence we can hope to get a contact sphere theorem using this strategy
for the pinching constant $1/4$.  Unfortunately, our tightness radius estimate
is not strong enough for this approach to work directly and extra topological
considerations were needed to prove Theorem \ref{CST-thm}. It may be that there
are better tightness radius estimates in pinched manifolds analogous to what
happens with the classical injectivity radius (there is no good general
injectivity radius estimate in odd dimensions and Klingenberg's estimate used
above really needs pinched curvature).

While our proof does not give any clues as to how one might improve the
pinching constant from $4/9$ to $1/4$, we notice that one might hope for an even
better pinching constant given that in dimension 3, Hamilton's Ricci flow
\cite{Hamilton82} allows one to prove that any pinching constant (or even
asking only the Ricci curvature to be positive) will suffice in the topological
sphere theorem.

\begin{question}
\label{CST}
What is the optimal curvature condition in the contact sphere theorem?
\end{question}

One might also ask for a higher dimensional analog of the contact sphere 
theorem, but at the moment we seem to know too little about contact topology
outside of dimension 3 for an approach similar to the one taken here to work (see \cite{EKM-higher} for some steps in this direction).

\section{Examples}\label{S:examples}
In this section we apply our theorems to several examples. In particular in the
first subsection we look at situations where Theorem~\ref{thm:hyp-criterion} can
be used to prove tightness of contact structures. In the following subsection we
consider overtwisted contact structures and compare the estimate we have on the
tightness radius to the actual tightness radius. 

\subsection{Examples in nonpositive curvature}
Our first goal is to provide examples when Theorem~\ref{thm:hyp-criterion} is
applicable. 

{\it Flat 3-torus:} We begin by investigating the well-known family of contact
structures $\xi_k$ on the 3-torus $T^3$ in standard coordinates, defined as
$\xi_k=\ker\alpha_k,$ for $k\in \mathbb{N}$, where 
$\alpha_k=\cos(k z)\,dx-\sin(k z)\,dy.$ The Reeb field is given by 
$R_{\alpha_k}=\cos(k z)\,\partial_x-\sin(k z)\,\partial_y.$
In the flat metric $g_{T^3}=dx^2+dy^2+dz^2$ on $T^3$, one may compute 
\[
 \ast\,d\alpha_k=k\,\alpha_k\,\quad \text{and} \quad \|\alpha_k\|=\|R_k\|=1.
\]
Therefore by definition the flat metric is compatible with all $\xi_k$. 
Thanks to Equation~\eqref{eq:m_g} we get $\mg=0$ and since
$\sec(g_{T^3})\equiv 0$, the Inequality~\eqref{eq:m_g-K} of Theorem
\ref{thm:hyp-criterion} holds and the theorem concludes that $\xi_k$ are
universally tight (as is well-known).

{\it Hyperbolic 3-space:} Consider the upper half space model of the 
hyperbolic 3-space i.e. $\mathbb{H}^3=\{(x,y,z)\in \R^3 | z>0\}$ equipped with
the metric
\[
 g_{\mathbb{H}^3}=\frac{1}{z^2}\bigl(dx^2+dy^2+dz^2\bigr).
\] 
Consider the same family of contact structures $\xi_k$ and contact forms
$\alpha_k$ as above. Since the hyperbolic metric $g_{\mathbb{H}^3}$ is conformal
to the flat metric, the Reeb field $R_k$ is orthogonal to $\xi_k$, thus
$g_{\mathbb{H}^3}$ is weakly compatible with each $\xi_k$. Let us determine the
``parameters'' of Theorem~\ref{thm:hyp-criterion}.  One may check that
\[
 \ast\,d\alpha_k= (kz)\, \alpha_k,
\]
and thus obtain that $ \rho  =\frac 1z$ and $\theta'=k\,z.$ Hence we see that
\[
\begin{split}
 d\ln(\rho) & =-\frac1z\,dz,\quad \nabla\ln\rho = -z\partial_z,\\
 d\ln(\theta') & = \frac1z dz\quad\text{and}\quad (\nabla\ln\theta')^\perp =0.
\end{split}
\]
Since $\sec(g_{\mathbb{H}^3})\equiv -1$ and  $
 \mg=\max_{\mathbb{H}^3}\|\nabla\ln\rho-(\nabla\ln\theta')^\perp\|=1,$
Inequality~\eqref{eq:m_g-K} of Theorem~\ref{thm:hyp-criterion} holds and the
theorem concludes again that $\xi_k$ are tight.

Curl eigenfields are of special interest in hydrodynamics and other branches of
physics \cite{EtnyreGhrist00}. As an example we point out is the the family of
1-forms:
\[
 \alpha'_k:=\cos(k\ln(z))\,dx-\sin(k\ln(z))\,dy,
\]
which are curl eigenfields on $\mathbb{H}^3$ with eigenvalues $\theta'=k$.

\medskip
{\it Product $\mathbb{R}\times\mathbb{H}^2$:} For the third nonpositive
curvature example let $M=\mathbb{R}\times\mathbb{H}^2$ be equipped
with the product metric 
\[
 g_{\mathbb{R}\times\mathbb{H}^2}=dt^2+\frac1y \bigl(dx^2+dy^2\bigr),
 \qquad (t,x,y)\in \mathbb{R}\times\mathbb{H}^2,\ y>0.
\]
\no We consider the following 1-form on $M$:
\[ 
 \beta=y^{\frac 12}\,dt+y^{-\frac12}\,dx,
\]
since $\beta$ is a positive rescaling of $dx+y\,dt$, $\xi=\ker\beta$ is just the
standard tight contact structure on $M\cong \R^3$. One may compute  
$\ast d\beta  =\frac12 \beta$
and that $ \rho =\|\beta\|^{-1}=(2\,y)^{-\frac12},$ and $\theta'=\frac 12.$
Therefore, $\beta$ defines a $\frac12$-curl eigenfield on $M$. Calculation of
the parameters of Theorem~\ref{thm:hyp-criterion} yields
\[
\begin{split}
 d\ln(\rho) & =-\frac{dy}{2\,y},
 \quad \nabla\ln\rho=-\frac12 y\partial_y,\quad  d\ln(\theta') = 0, \text{ and}\\
 \mg & =\max_M \|\nabla\ln\rho\|=\frac12.
\end{split}
\]
Since, $\sec(g_{\R\times\mathbb{H}^2})\leq 0$, the condition of
Theorem~\ref{thm:hyp-criterion} fails, and from Theorem
\ref{thm:weak-compatible} we can only get a tightness radius estimate
($\conv(g_{\R\times\mathbb{H}^2})=+\infty$) of
\[
 \tau(\R\times\mathbb{H}^2,\xi)\geq \text{ct}^{-1}_0\bigl(\frac 12\bigr)=2.
\]

\subsection{Overtwisted contact structures in flat 3-space.}\label{S:overt-num}
On $\R^3$ with the Euclidean metric $g=dr^2+r^2\,d\theta^2+dz^2$ in 
cylindrical coordinates $(r, \theta, z)$ we consider contact forms which are
invariant under vertical translations, rotations around the $z$-axis and
tangent to rays orthogonal to the $z$ axis:
\[
\alpha = a(r)\, dz + b(r)\, d\theta.
\]
\no This defines a smooth 1-form if $b(r) = r^2 \tilde b$ for some smooth
function $\tilde b$. Let $\vol$ be the usual euclidean volume form on $\R^3$ and
$\delta$ be defined by $\alpha \wedge d\alpha = \delta \, vol_g$. Here we will
assume that $\alpha$ is a positive contact form, i.e. $\delta > 0$. Observe that
\[
\delta = \frac{ab'-a'b}{r}, 
\quad
\|\alpha\| = \sqrt{a^2 + \left(\frac{b}{r}\right)^2},
\text{ and}\quad
R_\alpha = \frac{1}{\delta} \left(b'/r\partial_z - a'/r \partial_\theta\right).
\]

\no The line field normal to $\ker \alpha$ is directed by the vector field
$a \partial_z + \frac{b}{r^2} \partial_\theta$ dual to $\alpha$. 
So $\xi=\ker \alpha$ is weakly compatible with the Euclidean metric if
\begin{equation}
\label{eqn:sys}
\begin{cases}
	b' &= ra \\
	a' &= -\frac{b}{r}.
\end{cases}
\end{equation}
Note that the equation is sufficient but not necessary, we have considerably
more freedom: a priori we could have multiplied both right hand sides by any non
vanishing function of $r$.

A priori Equation \eqref{eqn:sys} is not well behaved at $r = 0$. But if a
solution exists then we have $b = -ra'$ so $b' = -ra'' - a'$ which when combined
with the first equation yields
\begin{equation}
\label{eqn:bessel}
ra'' + a' + ra = 0
\end{equation}
This is equivalent to $r^2a'' + ra' + r^2a = 0$ and $a$ continuous at zero.
We recognize the Bessel equation of order 0. Thus, for any choice of $b(0)$ we can
use the solution $a = a(0) J_0$, where $J_0$ is the zeroth Bessel function of
the first kind, and set $b = -ra'$ as required. Both $a$ and $b$ are analytic
functions defined for all $r \geq 0$. We may check that $a$ and $b$ cannot
both vanish simultaneously so that $\alpha$ is nonsingular, and a contact form. 
So we have a contact structure $\xi_\text{ot}=\ker\alpha$ which is overtwisted
at infinity because of the oscillatory behavior of $J_0(r),$ which is well known
to vanish for arbitrarily large $r.$ This must be compared to
Corollary~\ref{cor:ut} which implies that any contact structure compatible with
the flat metric on $\R^3$ is tight. Again we see that weak compatibility is much
more flexible.

\no Numerical experiments performed with the above solution suggest the
tightness radius estimate given by Theorem~\ref{thm:weak-compatible} is about
$0.15$ whereas we know the contact structure is tight when restricted to the
cylinder $D_r\times \R$ where $r$ is less than the first zero of $J'_0=-J_1.$
This first zero is slightly greater than 3.


\def\cprime{$'$} \def\cprime{$'$}

\end{document}